\documentclass[journal,twoside,web]{ieeecolor}
\usepackage{generic}
\usepackage{cite}
\usepackage{amsmath,amssymb,amsfonts,mathtools}
\usepackage{mathrsfs}
\usepackage{algorithmic}
\usepackage{graphicx}
\usepackage{algorithm,algorithmic}
\usepackage{textcomp}
\def\BibTeX{{\rm B\kern-.05em{\sc i\kern-.025em b}\kern-.08em
    T\kern-.1667em\lower.7ex\hbox{E}\kern-.125emX}}
\newtheorem{thm}{Theorem}

\newtheorem{lemma}[thm]{Lemma}

\newtheorem{prop}[thm]{Proposition}

\newtheorem{exmp}[thm]{Example}
\newtheorem{remark}[thm]{Remark}

\allowdisplaybreaks

\begin{document}
\title{Global Convergence of Oja's Component Flow for General Square Matrices and Its Applications}
\author{Daiki Tsuzuki and Kentaro Ohki, \IEEEmembership{Member, IEEE}
\thanks{This work was supported by JSPS KAKENHI Grant Numbers 23H01432 and 21K12097.}
\thanks{Daiki Tsuzuki is with the Graduate School of Informatics, Kyoto University, Kyoto, Japan (e-mail: tsuzuki.daiki.57a@bode.amp.i.kyoto-u.ac.jp)
}
\thanks{Kentaro Ohki is with the Graduate School of Informatics, Kyoto University, Kyoto, Japan (e-mail: ohki@i.kyoto-u.ac.jp).}
}

\maketitle

\begin{abstract}
In this study, the global convergence properties of the Oja flow, a continuous-time algorithm for principal component extraction, was established for general square matrices. The Oja flow is a matrix differential equation on the Stiefel manifold designed to extract a dominant subspace. Although its analysis has traditionally been restricted to symmetric positive-definite matrices, where it acts as a gradient flow, recent applications have extended its use to general matrices. In this non-symmetric case, the flow extracts the invariant subspace corresponding to the eigenvalues with the largest real parts. However, prior convergence results have been purely local, leaving the global behavior as an open problem. 
The findings of this study fill this gap by providing a comprehensive global convergence analysis, establishing that the flow converges exponentially for almost all initial conditions. 
We also propose a modification to the algorithm that enhances its numerical stability. As an application of this theory, we developed novel methods for model reduction of linear dynamical systems and the synthesis of low-rank stabilizing controllers. 
The study advances the theoretical understanding of the Oja flow and demonstrates its potential as a reliable and versatile tool for analyzing and controlling complex linear systems.
\end{abstract}


\section{Introduction}
\label{sec:introduction}

\IEEEPARstart{P}{rincipal} component analysis (PCA) is a foundational technique for dimensionality reduction and feature extraction, with broad applications across technology and science \cite{moore1981principal,bishop2006pattern,jolliffe2016principal}. For large-scale and streaming datasets, online algorithms for PCA are essential. A prominent class of such algorithms is principal and minor component flows, which are recursive methods that have gained significant attention in statistics, machine learning \cite{Plaut2018,cunningham2022principal,wright2022high,wainwright2019high}, signal processing \cite{tanaka2005generalized,peng2006convergence}, and control theory \cite{hasan2006stability,manton2005dual,borkar2012oja,benner2015survey,bonnabel2012geometry,yamada2020lowrank,yamada2021comparison,TsuzukiOhki2024,TsuzukiOhki2024a}. Although these algorithms are valued for their computational efficiency, the theoretical conditions guaranteeing their convergence, particularly for general non-symmetric matrices, remain incompletely understood.

\subsection{Background of the Oja Flow}
A canonical example of a principal component flow is the {\em Oja flow} \cite{oja1982simplified,oja1985stochastic,oja1989neural}:
\begin{align}
    \varepsilon \frac{d}{dt}U(t) = (I_n - U(t)U(t)^{\top}) AU(t), \ U(0) \in \mathrm{St}(r,n),
    \label{eq:Oja_flow}
\end{align}
where $r\leq n$, $A \in \mathbb{R}^{n\times n}$, $\varepsilon \in (0,1]$ is a rate-controlling parameter, $I_n$ is the $n \times n$ identity matrix, and $\mathrm{St}(r,n) := \{ X \in \mathbb{R}^{n\times r} \mid X^{\top}X = I_r\}$ is the Stiefel manifold. The corresponding minor component flow, which seeks the least dominant subspace, is obtained by replacing $A$ with $-A$.

Although the Oja flow is widely used in machine learning to extract the dominant singular subspaces of symmetric positive-definite matrices, a case supported by strong theoretical guarantees \cite{chen1998global,chen2001unified,yoshizawa2001convergence,manton2005dual}, its analysis and application for general square matrices have been limited \cite{kuo2025asymptotic}.

This paper provides a comprehensive convergence analysis of the Oja flow for general matrices and explores its applications in control theory. Our previous work \cite{TsuzukiOhki2024} showed that the Oja flow can extract the invariant subspace corresponding to the $r$ eigenvalues with the largest real parts. However, that analysis was restricted to local convergence, and the estimate of the domain of attraction was conservative. This work extends these results to establish global convergence properties.

\subsection{Related Work}
Prior studies on the Oja flow have primarily focused on symmetric positive-definite matrices, yielding several key results:
\begin{enumerate}
    \item The existence and uniqueness of the solution to the \eqref{eq:Oja_flow} were established in \cite{yan1994global}.
    \item For any full-rank initial matrix $U(0) \in \mathbb{R}^{n\times r}$, the solution $U(t)$ converges to the Stiefel manifold $\mathrm{St}(r,n)$ \cite{yan1994global}.
    \item For almost all initial conditions, the solution converges to the subspace spanned by the eigenvectors corresponding to the $r$ largest eigenvalues, provided the eigenvalues are distinct \cite{chen1998global}.
\end{enumerate}
These results extend to any symmetric matrix by a spectral shift, $A \to A + aI_n$, $a>0$, which renders the matrix positive-definite without altering the eigenvectors. Extensions have broadly followed two directions: modifying the Oja flow itself or analyzing the original flow for a more general class of matrices. This paper pursues the latter direction, following preliminary work in \cite{hasan2006stability,bonnabel2012geometry}.

Concurrently with our work, Kuo et al. \cite{kuo2025asymptotic} have also established the existence and uniqueness of the solution for general square matrices and analyzed its exponential convergence on the Stiefel manifold under a certain spectral gap condition. Although this significant work addresses convergence on the Stiefel manifold, a rigorous analysis of the global domain of attraction and the solution behavior for initial conditions outside the manifold remains an open area. The aim of our study also included addressing these specific gaps.

\subsection{Contributions}
The main contributions of this paper are as follows:
\begin{enumerate}
    \item A comprehensive convergence analysis of the Oja flow for general real square matrices (Theorems \ref{prop:positive_part_convergence} and \ref{thm:convergence_Oja_flow}). We prove exponential convergence to the dominant invariant subspace under a mild eigenvalue separation condition and propose a modification that ensures numerical stability by guaranteeing convergence to the Stiefel manifold.
    
    \item Estimation of the domain of attraction (Theorem \ref{thm:volume_of_basin}). We show that the domain of attraction encompasses almost the entire manifold, indicating that the flow converges for almost all initial values.
    
    \item Applications to control theory. We developed a novel framework for model reduction and low-rank controller synthesis. We show that the proposed reduction method preserves key system properties (observability and controllability) and can be used to design stabilizing controllers for large-scale systems with low-dimensional unstable manifolds.
\end{enumerate}
Although we focus on real matrices, our results extend directly to complex matrices.

\subsection{Organization and Notation}
This paper is organized as follows. Section \ref{sec:summary} reviews existing results. Section \ref{sec:main_Oja_flow} presents our main theoretical contributions on the convergence of the Oja flow, including its numerical stability and domain of attraction. Section \ref{sec:applications} demonstrates applications to model reduction and controller synthesis. Key theoretical results are followed by numerical examples for validation.

\subsubsection*{Notation}
The sets of real and complex numbers are $\mathbb{R}$ and $\mathbb{C}$. The set of $n \times m$ real matrices is $\mathbb{R}^{n\times m}$. $I_n$ is the $n \times n$ identity matrix, and $0_{n,m}$ is the $n \times m$ zero matrix. For simplicity, we denote $0_{n} = 0_{n,n}$; if the dimension is trivial, $0$ is also used. 
$A^{\top}$ and $A^{\dagger}$ denote the transpose and Hermitian conjugate of a matrix $A$. $\|x\|$ is the Euclidean norm for a vector $x$ and $\| A \| _{\rm ind}$ is its induced norm for a matrix $A$. For a symmetric matrix $A$, $A > 0$ ($A \geq 0$) indicates it is positive-definite (semidefinite). $A^{1/2}$ denotes the unique positive-semidefinite square root. The eigenvalues of a square matrix $A \in \mathbb{C}^{n \times n}$ are ordered such that $\mathrm{Re}(\lambda_1(A)) \geq \dots \geq \mathrm{Re}(\lambda_n(A))$. The corresponding (generalized) eigenvector is $\psi_i(A)$, and the matrix of eigenvectors is $\Psi(A) := [\psi _{1}(A) , \dots , \psi _{n}(A)]$. The Stiefel manifold is $\mathrm{St}(r,n) := \{X \in \mathbb{R}^{n \times r} \mid X^{\top}X = I_r\}$.

\section{Summary of Existing Results}
\label{sec:summary}

Existing theoretical work on the key properties of the Oja flow \eqref{eq:Oja_flow} can be broadly divided according to whether the initial value $U(0)$ lies on the Stiefel manifold $\mathrm{St}(r,n)$.

\subsection{Convergence Analysis in Euclidean Space}
	
The first category of results concerns the solution behavior of the Oja flow \eqref{eq:Oja_flow} in Euclidean space. 
It is well established that if the initial condition $U(0)$ is on the Stiefel manifold, i.e., $U(0) \in \mathrm{St}(r,n)$, then the solution $U(t)$ remains on $\mathrm{St}(r,n)$ for all $t \geq 0$ and any matrix $A \in \mathbb{R}^{n\times n}$ \cite[\S II.B]{TsuzukiOhki2024}. 
In practice, however, numerical errors can cause $U(t)$ to deviate from the manifold. 
If $\mathrm{St}(r,n)$ were a stable invariant set, a simple numerical integrator such as the forward Euler scheme would keep the solution in a neighborhood of the manifold. However, the stability conditions for general matrices have not been fully investigated. 
For the specific case of symmetric positive-semidefinite matrices, the following result was shown in \cite{yan1994global}:

\begin{prop}[{\cite[Prop. 3.1]{yan1994global}}]
\label{prop:yan1994global_prop3.1}
Let $U(0) \in \mathbb{R}^{n\times r}$ be of full rank.   
Then, for any symmetric positive-definite matrix $A \in \mathbb{R}^{n\times n}$, the solution $U(t)$ of \eqref{eq:Oja_flow} converges exponentially to $\mathrm{St}(r,n)$ as $t\to \infty$. 
\end{prop}
	
This result implies that for positive-definite matrices $A$, $\mathrm{St}(r,n)$ is a stable invariant set for the Oja flow.
In contrast to Proposition \ref{prop:yan1994global_prop3.1}, if $A$ is not positive-definite, the Oja flow may fail to converge to $\mathrm{St}(r,n)$ for an initial condition $U(0) \notin \mathrm{St}(r,n)$.
The following example illustrates this idea.

\begin{exmp}\label{exam:negative_eigenvalue}
Consider a system with
\begin{align*}
	A = \mathrm{diag}(1 ,\ -1) ,\quad 
	U(0) = \begin{bmatrix} 0 & a \end{bmatrix} ^{\top}, \ a>1. 
\end{align*}
The initial derivative for the Oja flow \eqref{eq:Oja_flow} is
\begin{align*}
	\frac{d}{dt} U(t) \bigg|_{t=0} = \frac{1}{\varepsilon} \begin{bmatrix} 1 & 0 \\ 0 & a^2 -1 \end{bmatrix} U(0) 
	= 
	\frac{1}{\varepsilon} 
	\begin{bmatrix} 0 \\ a (a^{2}-1) \end{bmatrix},
\end{align*}
which indicates that the second component of the solution $U(t)$ monotonically increases from its initial value.
Hence, $U(t)$ does not converge to $\mathrm{St}(r,n)$. 
\end{exmp}

Chen et al. \cite{chen1998global} extended Proposition \ref{prop:yan1994global_prop3.1} to positive-semidefinite matrices and rank-deficient initial conditions, showing that the Oja flow converges to a point in $\mathbb{R}^{n\times r}$ that depends on $U(0)$.
Hasan \cite[Variation 4]{hasan2006stability} claimed that convergence to $\mathrm{St}(r,n)$ could be ensured for any square matrix $A$ and any non-zero initial value $U(0)$ by replacing $A$ with $A+aI_{n}$ for a sufficiently large constant $a>0$. 
This modification does not alter the vector field of the Oja flow \eqref{eq:Oja_flow} when restricted to $U(t) \in \mathrm{St}(r,n)$.
However, a non-zero initial value is not a sufficient condition for convergence.
For instance, if the columns of $U(0) \in \mathbb{R}^{n\times 2}$ are identical, i.e., $U(0) = x_0 \begin{bmatrix} 1 & 1 \end{bmatrix}$, then the solution retains this structure, $U(t) = x(t) \begin{bmatrix} 1 & 1 \end{bmatrix}$, where $x(t)$ evolves according to
\begin{align*}
	\frac{d}{dt}x (t) = (I_{n} - 2 x(t) x(t)^{\top} ) A x(t) ,\ x(0) =x_{0}. 
\end{align*}
Therefore, $U(t)$ preserves its initial rank and cannot converge to $\mathrm{St}(r,n)$. 
The strategy of shifting the matrix spectrum is compelling, but further analysis is needed because $\mathrm{St}(r,n)$ is generally not a stable invariant set for non-positive-definite matrices.
As demonstrated in Example \ref{exam:negative_eigenvalue} and discussed in \cite{chen2001unified}, if $A$ is negative-definite, solutions initiating outside the Stiefel manifold may diverge.
The Stiefel manifold is a generalization of a sphere, and this potential instability implies that numerical implementations of the Oja flow require periodic normalization to mitigate the accumulation of numerical errors.
To address this limitation, Chen et al. \cite{chen1998unified, chen2001unified} proposed an alternative algorithm for principal and minor component extraction that does not require normalization.
In this study, however, we focused on the Oja flow \eqref{eq:Oja_flow} because of its lower computational complexity compared to the algorithm in \cite{chen1998unified}.

	The stability properties of \eqref{eq:Oja_flow} in $\mathbb{R}^{n \times r}$, particularly for general non-symmetric matrices, remain incompletely understood.
The aim of this study was to fill this gap by providing a rigorous analysis of the invariance and stability of $\mathrm{St}(r,n)$ for any $A \in \mathbb{R}^{n\times n}$, as described in Section \ref{sec:Oja_flow_Euclidean}.

\subsection{Convergence Analysis on the Stiefel Manifold}

As a principal component extraction algorithm, the Oja flow is expected to converge to the subspace spanned by the eigenvectors corresponding to the dominant eigenvalues.
This property is known to hold for positive-definite matrices. 

\begin{prop}[\cite{yan1994global}]\

\begin{enumerate}
    \item (\cite[Theorem 5.1]{yan1994global}) If $A$ is symmetric positive-definite with $\lambda_{r} > \lambda_{r+1}$, then for any initial condition $U(0)$ satisfying $\mathrm{det}( U(0)^{\top} \Psi_{r} ) \neq 0$, the solution of the Oja flow converges to the set
    \begin{align}
        \mathcal{U}_{r}^{\prime} := \left\{ \Psi(A) \begin{bmatrix}
        K_{r} \\ 0_{n-r,r}
        \end{bmatrix}
        \in \mathrm{St}(r,n) \ | \ K_{r}\in \mathbb{R}^{r\times r}
        \right\},
        \label{eq:StableEquilibriumSets}
    \end{align}
    where $\Psi(A)$ is the matrix of eigenvectors of $A$, and $\Psi_{r} := [\psi_{1}(A) ,\dots , \psi_{r}(A)]$.
    
    \item (\cite[Corollary 5.1]{yan1994global}) The initial condition $\mathrm{det}( U(0)^{\top} \Psi_{r} ) \neq 0$ holds for almost all $U(0) \in \mathrm{St}(r,n)$.
\end{enumerate}
\end{prop}
Notably, these results also apply for symmetric matrices through the use of a spectral shift. 
Recently, these results were extended to general matrices $A \in \mathbb{R}^{n\times n}$.

\begin{prop}[\cite{TsuzukiOhki2024,TsuzukiOhki2024a,kuo2025asymptotic}] \label{prop:TsuzukiOhki2024}
Assume that the eigenvalues of $A$ are ordered such that $\mathrm{Re}(\lambda_{r}) > \mathrm{Re}(\lambda_{r+1})$. Then, the following properties hold:
\begin{enumerate}

    \item (\cite[Prop. 2]{TsuzukiOhki2024}) The equilibrium sets of \eqref{eq:Oja_flow} are given by
    \begin{align}
        \mathcal{U}_{r,P} := \left\{ \Psi(A) P \begin{bmatrix}
        K_{r} \\ 0_{n-r,r}
        \end{bmatrix}
        \in \mathrm{St}(r,n) \ | \ K_{r}\in \mathbb{C}^{r\times r}
        \right\} ,
        \label{eq:EquilibriumSets}
    \end{align}
    where $P\in \mathbb{R}^{n\times n}$ is a permutation matrix such that $P^{\top} \Lambda P$ is a block-diagonal with blocks of size $r\times r$ and $(n-r)\times (n-r)$. For simplicity, we denote $\mathcal{U}_{r} = \mathcal{U}_{r,I_{n}}$. Note that for any $\bar{U} \in \mathcal{U}_{r,P}$, the matrix $\bar{U}W$ is also in $\mathcal{U}_{r,P}$ for any orthogonal matrix $W \in \mathbb{R}^{r\times r}$.

    \item (\cite[Prop. 3]{TsuzukiOhki2024}) For any $\bar{U} \in \mathcal{U}_{r,P}$, the eigenvalues of $\bar{U}^{\top} A \bar{U}$ are a permutation of the first $r$ eigenvalues of $A$, i.e., $\{ \lambda_{i}(\bar{U}^{\top} A \bar{U}) \}_{i=1}^{r} = \{ \lambda_{\mathcal{I}_{P}(i)} (A) \}_{i=1}^{r}$, where $\mathcal{I}_{P}$ is the permutation associated with $P$. In particular, for $\bar{U} \in \mathcal{U}_{r}$, we have $\lambda_{i}(\bar{U}^{\top}A\bar{U}) = \lambda_{i}(A)$ for $i= 1,\dots ,r$.

    \item (\cite[Lemma 1]{TsuzukiOhki2024a}) For any $\bar{U} \in \mathcal{U}_{r,P}$ and any integer $p \ge 0$, the subspace spanned by the columns of $\bar{U}$ is an invariant subspace of $A$, satisfying $A^p \bar{U} = \bar{U} (\bar{U}^{\top}A\bar{U})^{p}$. This result implies $\mathrm{e}^{A} \bar{U} = \bar{U} \mathrm{e}^{\bar{U}^{\top}A \bar{U}}$.

    \item (\cite[Thm. 1]{TsuzukiOhki2024}) The set $\mathcal{U}_{r}$ is the unique asymptotically stable equilibrium set, which means any trajectory starting sufficiently close to $\mathcal{U}_{r}$ remains in its neighborhood and converges to $\mathcal{U}_{r}$.
    
    \item (\cite{kuo2025asymptotic}) 
    For any initial condition $U_{0}$ in the set
     \begin{align}
        \label{eq:Kuo2025asymptotic}
        \mathcal{V}_{r} := \Bigg{\{} 
        & \Psi(A) \begin{bmatrix}
        K_{r} \\ K_{\perp}
        \end{bmatrix}
        \in \mathrm{St}(r,n)
        \ | \ K_{r} \in \mathbb{C}^{r\times r} , \nonumber 
        \\ &  
        K_{\perp} \in \mathbb{C}^{(n-r)\times r}, \ 
        \mathrm{rank}(K_{r}) =r 
        \Bigg{\}},
    \end{align}
    the solution $U(t)$ of \eqref{eq:Oja_flow} converges exponentially to $\mathcal{U}_{r}$.

\end{enumerate}
\end{prop}

However, these properties are insufficient for a complete convergence analysis for two main reasons:
\begin{enumerate}
	\item 
	In many practical applications, the strict spectral gap condition at the desired subspace dimension, $\mathrm{Re}(\lambda_r) > \mathrm{Re}(\lambda_{r+1})$, cannot be guaranteed or is difficult to verify. A convergence analysis under a milder assumption is therefore required.

    \item The existence of other invariant sets, such as limit cycles, has not been ruled out \cite{TsuzukiOhki2024}.

\end{enumerate}

The milder assumption was addressed in this study by investigation of the case where a spectral gap exists at an index $m \geq r$, such that $\mathrm{Re}(\lambda_m) > \mathrm{Re}(\lambda_{m+1})$.  
Furthermore, the analysis of the domain of attraction shows that no other meaningful invariant set exists.

\begin{figure}[!htbp]
    \centering
    \includegraphics[keepaspectratio,width=\linewidth]{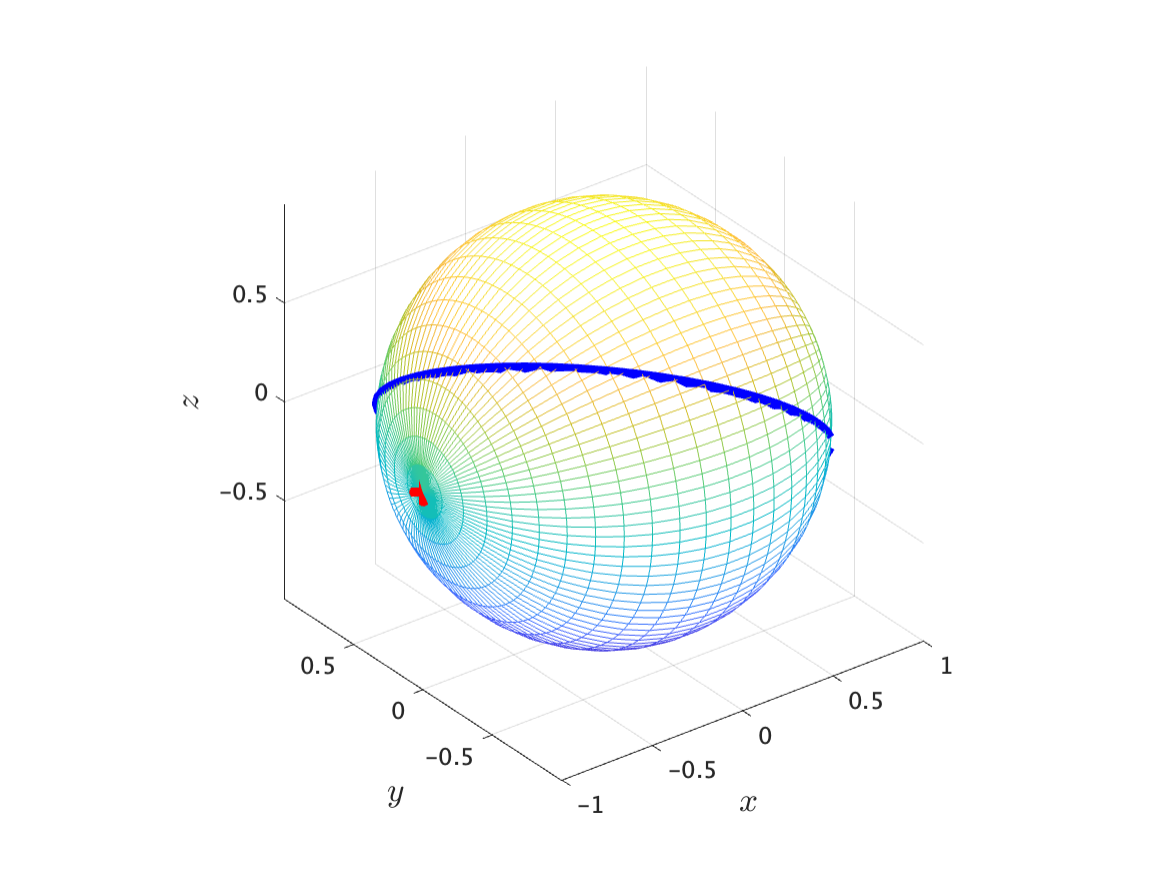} 
    \caption{Sphere representing $\mathrm{St}(1,3)$.  The red markers represent $\mathcal{U}_{1}=\{ \pm \psi _{1}\}$ and the blue line represents the unit circle $\mathrm{St}(1,3) \cap \mathrm{span}\{ \psi _{2}, \psi _{3} \}$ in Example \ref{exmp:St13}.  }
    \label{fig:sphere_case}
\end{figure}

\begin{exmp}[Visualization of the set $\mathcal{V}_{r}$] \label{exmp:St13}
Consider the Oja flow with the matrix:
\begin{align*}
    A = \begin{bmatrix}
    1 & 1 & 2 \\
    0 & 0 & 1 \\
    0 & 0 & -1
    \end{bmatrix}.
\end{align*}
The eigenvectors of $A$ are
\begin{align*}
    \psi_{1} = \begin{bmatrix} 1 \\ 0 \\ 0 \end{bmatrix}, \quad
    \psi_{2} = \frac{1}{\sqrt{2}}\begin{bmatrix} 1 \\ -1 \\ 0 \end{bmatrix}, \quad
    \psi_{3} = \frac{1}{3} \begin{bmatrix} 1 \\ 2 \\ -2 \end{bmatrix}.
\end{align*}
These vectors are not mutually orthogonal. For $(n,r) = (3,1)$, the stable equilibrium set is $\mathcal{U}_{1} = \{ \pm \psi_{1}\}$, whereas the unstable equilibria are $\{ \pm \psi_{2}\} $ and $\{ \pm \psi_{3}\}$.
Fig. \ref{fig:sphere_case} depicts the stable equilibria $\mathcal{U}_{1}$ (red markers) and the unit circle formed by the intersection of $\mathrm{St}(1,3)$ and the plane $\mathrm{span}\{ \psi_{2}, \psi_{3} \}$ (blue line). 
The initial set \eqref{eq:Kuo2025asymptotic} in Prop. \ref{prop:TsuzukiOhki2024} is then the unit sphere without the blue line, which means that almost all initial values converge to $\mathcal{U}_{1}$. 
This geometric observation suggests that the domain of attraction has full measure on the manifold, a claim we formalized as discussed in Section \ref{sec:basin}.

\end{exmp}

\subsection{Other Related Work}
		
		Several variants of the Oja flow have been proposed \cite{mathieu2020riemannian,cunningham2022principal,bonnabel2024low,kuo2023orthogonal}. 
A key property of the standard Oja flow \eqref{eq:Oja_flow} is its invariance under right-multiplication by an orthogonal matrix. 
Specifically, if $U(t)$ is a solution, then for any differentiable orthogonal matrix $W(t) \in \mathbb{R}^{r \times r}$, the trajectory $U'(t) := U(t)W(t)$ spans the same subspace as $U(t)$ at every instant, because $U'(t)U'(t)^\top = U(t)W(t)W(t)^\top U(t)^\top = U(t)U(t)^\top$. 
Consequently, $U'(t)$ converges to the same invariant subspace as $U(t)$.
The dynamics of $U'(t)$ can be expressed by augmenting the Oja flow with a term related to the derivative of $W(t)$. Because $W(t)$ is orthogonal, its derivative satisfies $\varepsilon \frac{d}{dt}W(t) = W(t)S(t)$ for some skew-symmetric matrix $S(t) = -S(t)^{\top} \in \mathbb{R}^{r\times r}$. 
The resulting dynamics for $U'(t)$ are given by
\begin{align*}
	\varepsilon \frac{d}{dt}U'(t) = (I_{n} - U'(t) U'(t)^{\top} ) AU'(t) + U'(t) S(t).
\end{align*}
This formulation enables modifications that impose additional structure on the solution without altering the fundamental subspace dynamics.

An important application of this principle is the continuous-time reduced QR algorithm \cite{frank2018detectability,tranninger2019detectability}. 
Setting $\varepsilon = 1$ for simplicity, this algorithm is described by the following coupled differential equations:
\begin{align}
	\frac{d}{dt} U_{\rm qr}(t) &= (I_{n} - U_{\rm qr}(t) U_{\rm qr}(t)^{\top}) A U_{\rm qr}(t) + U_{\rm qr}(t) S(t), \label{eq:QR_Oja_flow} \\
	\frac{d}{dt} R(t) &= B(t) R(t), \label{eq:reducedQR_R}
\end{align}
with initial conditions $U_{\rm qr}(0) \in \mathrm{St}(r,n)$ and $R(0) \in \mathbb{R}^{r\times r}$, where $R(0)$ is an upper triangular matrix. 
The skew-symmetric matrix $S(t)$ is chosen at each instant such that the matrix $B(t) := U_{\rm qr}(t)^{\top} A U_{\rm qr}(t) - S(t)$ is forced to be upper triangular.
If $A$ has real eigenvalues, the diagonal elements of the solution $R(t)$ converge to the $r$ eigenvalues of $A$ with the largest real parts.

This property makes the reduced QR algorithm a widely used tool for estimating Lyapunov exponents and dynamical system spectra, in a manner analogous to the full continuous-time QR algorithm (i.e., the case where $r=n$) \cite{dieci1997compuation,dieci2007lyapunov,bridges2001computing}.
The algorithm defined by \eqref{eq:QR_Oja_flow} and \eqref{eq:reducedQR_R} has found applications in Kalman--Bucy filtering for linear time-varying systems and in nonlinear observers \cite{tranninger2022detectability}.
Although this algorithm is designed for general square matrices and extensive numerical evidence suggests that it effectively extracts the dominant subspace, a rigorous theoretical analysis of its convergence properties is still lacking, even for the time-invariant case.

\section{Convergence Analysis for Oja Flow}
\label{sec:main_Oja_flow}

    This section provides the theoretical results of the Oja flow \eqref{eq:Oja_flow}.  
    First, we demonstrate how to stabilize the Stiefel manifold $\mathrm{St}(r,n)$ and establish the convergence rate in Euclidean space $\mathbb{R}^{n\times r}$. 
    Next, we clarify the convergence on $\mathrm{St}(r,n)$, and subsequently establish the domain of attraction.  
    Finally, we provide guidance on how to efficiently increase or decrease $r$.  

\subsection{Convergence to the Stiefel Manifold}
\label{sec:Oja_flow_Euclidean}

	In this section, we extend the convergence results for the Oja flow from symmetric positive-definite matrices, as presented in Proposition \ref{prop:yan1994global_prop3.1}, to the general case of arbitrary square matrices.
First, we present a lemma that generalizes Lemma 2.2 of \cite{yan1994global}.
Although Hasan previously established the asymptotic convergence of \eqref{eq:Oja_flow} for general square matrices \cite[Variation 4]{hasan2006stability}, our analysis results demonstrated exponential convergence.
The following result from the literature is instrumental to our proof.

\begin{lemma}[{\cite{sasagawa1982finite}}]
\label{lemma:sasagawa1982finite}
Consider the Riccati differential equation:
\begin{align*}
	\frac{d}{dt}P(t) = AP(t) + P(t) A^{\top} - P(t) Q P(t) + R,
\end{align*}
where $A, Q=Q^{\top}, R=R^{\top} \in \mathbb{R}^{n\times n}$, with $Q \geq 0$, $R\geq 0$, and $P(0)=P_{0} \geq 0$.
Let $\bar{P} \in \mathbb{R}^{n\times n}$ be an equilibrium solution.
Then, the solution is given by $P(t) = Y(t)X(t)^{-1}$ for all $t \in [0,t_{\max})$, where
\begin{align*}
	X(t) = &\mathrm{e}^{-\tilde{A}^{\top}t} 
	\left(
	I_{n} + \int _{0}^{t} \mathrm{e}^{\tilde{A}^{\top}s}  R  \mathrm{e}^{\tilde{A}s} ds (P_{0} - \bar{P})
	\right),
	\\
	Y(t) =& \bar{P} \mathrm{e}^{-\tilde{A}^{\top}t} 
	\nonumber \\ & 
	+ \left(
	\bar{P} \int _{0}^{t} \mathrm{e}^{-\tilde{A}^{\top}(t-s)}  R  \mathrm{e}^{\tilde{A}s} ds + \mathrm{e}^{\tilde{A}t}
	\right) (P_{0} - \bar{P}),
\end{align*}
with $\tilde{A} := A - \bar{P}R$, and $t_{\max} := \inf \{ t \geq 0 \ | \ \mathrm{det}(X(t)) = 0 \}$.
\end{lemma}

Using Lemma \ref{lemma:sasagawa1982finite}, we first establish that the rank of the solution matrix is preserved.

\begin{lemma}
\label{lemma:rank_U(t)}
Assume that the symmetric part of $A \in \mathbb{R}^{n\times n}$ is positive-definite, i.e., $A_{\rm sym} > 0$. Let $U(0) \in \mathbb{R}^{n\times r}$ be of full rank.
Then, the solution $U(t)$ of \eqref{eq:Oja_flow} remains full rank for all $t\geq 0$.
\end{lemma}

\begin{proof}
Let $P(t) := U(t)U(t)^{\top}$. Without loss of generality, we set $\varepsilon = 1$. The evolution of $P(t)$ is governed by the Riccati differential equation:
\begin{align*}
	\frac{d}{dt}P(t) = A P(t) + P(t) A^{\top} - 2 P(t) A_{\rm sym} P(t),
\end{align*}
with the initial condition $P(0) = P_{0} := U(0)U(0)^{\top}$.
Notably, $\bar{P} = 0_{n,n}$ is an equilibrium solution. Applying Lemma \ref{lemma:sasagawa1982finite} with $\bar{P}=0_{n}$, $Q=2A_{\rm sym}$, and $R=0_{n}$ yields the solution components
\begin{align*}
	X(t) =& \mathrm{e}^{-A^{\top} t} \left( I_{n} + G(t) P_{0} \right), \quad 
	Y(t) = \mathrm{e}^{At} P_{0},
\end{align*}
where $G(t) := 2\int _{0}^{t} \mathrm{e}^{A^{\top}s} A_{\rm sym} \mathrm{e}^{As} ds$.
Because $A_{\rm sym} > 0$, the integrand is positive-definite, and thus $G(t)>0$ for any $t>0$.
The eigenvalues of $G(t)P_{0}$ are the same as those of the symmetric matrix $G(t)^{1/2} P_{0} G(t)^{1/2}$, which are non-negative.
This result implies that all eigenvalues of $I_n + G(t)P_0$ are greater than or equal to 1 for all $t\geq 0$. Therefore, $X(t)$ is invertible for all $t\geq 0$, and the solution is
\begin{align}
	P(t) = Y(t) X(t)^{-1} = \mathrm{e}^{At} P_{0} \left( I_{n} + G(t) P_{0} \right)^{-1} \mathrm{e}^{A^{\top} t}.
	\label{eq:explicit_solution_Riccati}
\end{align}
Because $\mathrm{e}^{At}$ is always invertible, the rank of $P(t)$ is equal to the rank of $P_0$. This result implies that the rank of $U(t)$ is preserved for all $t\geq 0$.
\end{proof}

The following lemma, a generalization of Theorem 2.2 in \cite{yan1994global}, characterizes the evolution of the singular values of $U(t)$.

\begin{lemma}\label{lemma:svd_oja_solution}
Assume that $A_{\rm sym} > 0$. Let $U(0) \in \mathbb{R}^{n\times r}$ have full rank, and let $\sigma_{1}(t) \geq \dots \geq \sigma_{r}(t)$ denote the singular values of $U(t)$.
Then, the following hold:
\begin{enumerate}
    \item $\lim_{t\to \infty} \sigma_{i}(t) = 1$ for all $i=1,\dots,r$.
    \item If $\sigma_{i}(0) \geq 1$ for all $i$, then each $\sigma_{i}(t)$ is a non-increasing function of $t$.
    \item If $\sigma_{i}(0) \leq 1$ for all $i$, then each $\sigma_{i}(t)$ is a non-decreasing function of $t$.
    \item $\sigma_{r}(t) \geq \alpha$ for all $t\geq 0$, where $\alpha := \min\{1, \sigma_{r}(0)\}$.
\end{enumerate}
\end{lemma}

\begin{proof}
Let $P(t) = U(t)U(t)^{\top}$, with its solution given by \eqref{eq:explicit_solution_Riccati}. Without loss of generality, set $\varepsilon = 1$.
The matrix $G(t)$ in \eqref{eq:explicit_solution_Riccati} can be rewritten as
$G(t) = \int_{0}^{t} \frac{d}{ds}(\mathrm{e}^{A^{\top}s} \mathrm{e}^{As}) ds = H(t) - I_{n}$,
where $H(t) := \mathrm{e}^{A^{\top}t} \mathrm{e}^{At}$.
Take an orthogonal matrix $Q \in \mathbb{R}^{n\times n}$ such that
\begin{align*}
Q U(0) U(0)^{\top} Q^{\top} =& \mbox{blk-diag} (L,0_{n-r}),\\
L :=& \mathrm{diag}\begin{pmatrix} \sigma _{1}(0)^2 ,& \dots , & \sigma _{r}(0)^2 \end{pmatrix}.
\end{align*}
Then, $P(t)$ takes the form
\begin{align*}
P(t) =&\mathrm{e}^{At} Q^{\top} \mbox{blk-diag} (L,0_{n-r})
\\
& \times Q \{ I_{n} + (H(t) - I_{n}) P_{0} \} ^{-1} Q^{\top} Q \mathrm{e}^{A^{\top}t}
\\
=&\mathrm{e}^{At} Q^{\top} \mbox{blk-diag} (L,0_{n-r})
\\
& \times \{ I_{n} + (\Sigma(t) - I_{n})\mbox{blk-diag} (L,0_{n-r}) \} ^{-1} Q \mathrm{e}^{A^{\top}t}
\\
=& \mathrm{e}^{At} Q^{\top}
\begin{bmatrix}
\{ L^{-1} - I_{r} + \Sigma _{11}(t) \} ^{-1} & 0
\\
0 & 0
\end{bmatrix}
Q \mathrm{e}^{A^{\top}t},
\end{align*}
where $\Sigma (t) := QH(t)Q^{\top} $ and $\Sigma _{11}(t) \in \mathbb{R}^{r\times r}$ is the top-left block of $\Sigma(t)$.
Because $P(t)$ is symmetric positive semidefinite, its nonzero singular values are the square roots of its nonzero eigenvalues. These eigenvalues are the same as those of the following similarity transformed matrix:
\begin{align*}
&Q^{\top} \mathrm{e}^{A^{\top}t} \mathrm{e}^{At} Q
\begin{bmatrix}
\{ L^{-1} - I_{r} + \Sigma _{11}(t) \} ^{-1} & 0
\\
0 & 0
\end{bmatrix}
\\
=&
\begin{bmatrix}
\Sigma _{11}(t) \{ L^{-1} - I_{r} + \Sigma _{11}(t) \} ^{-1} & 0
\\
\Sigma _{21}(t) \{ L^{-1} - I_{r} + \Sigma _{11}(t) \} ^{-1} & 0
\end{bmatrix}
,
\end{align*}
where $\Sigma _{21}(t) \in \mathbb{R}^{(n-r) \times r}$ is the bottom-left block of $\Sigma(t)$.
Because $\Sigma (t) >0$, $\Sigma _{11}(t) >0$. Thus, the nonzero eigenvalues of $P(t)$ are the eigenvalues of
\begin{align*}
& \Sigma _{11}(t) \{ L^{-1} - I_{r} + \Sigma _{11}(t) \} ^{-1}
\\
=& \{ ( L^{-1} - I_{r} ) \Sigma _{11}(t)^{-1} + I_{r} \} ^{-1}.
\end{align*}
Recall that $\mathrm{Re} (\lambda _{i}(A)) \geq \lambda _{n}(A_{\mathrm{sym}}) > 0$ for all $i$ from the assumption $A_{\mathrm{sym}} >0$. Therefore, $H(t)$ is monotonically increasing in the Loewner order and will diverge as $t \to \infty$. 
This result implies that $\Sigma _{11}(t)$ also diverges, making $\Sigma _{11}(t) ^{-1} $ approach the zero matrix as $t\to \infty$. This result immediately shows that all the nonzero singular values of $U(t)$ approach 1, and hence, statement 1 holds true.

The monotonicity of $\Sigma _{11}(t)$ is then used to establish statements 2, 3, and 4 by following the same arguments as in the proof of Theorem 2.2 of \cite{yan1994global}.
\end{proof}
	
The condition $A_{\rm sym} > 0$ may seem restrictive, but as shown in the next theorem, it can always be satisfied by a simple modification of the Oja flow \eqref{eq:Oja_flow} that does not alter its dynamics on the Stiefel manifold.

\begin{thm} \label{prop:positive_part_convergence}
For a given $A \in \mathbb{R}^{n\times n}$, choose a scalar $a \geq 0$ such that $A_{\rm sym} + a I_{n}$ is positive-definite.
Let the initial condition $U(0) \in \mathbb{R}^{n\times r}$ be a full-rank matrix.
Then, the solution of the modified Oja flow
\begin{align*}
    \varepsilon \frac{d}{dt} U(t) = (I_{n} - U(t) U(t)^{\top}) (A+aI_n) U(t)
\end{align*}
converges exponentially to $\mathrm{St}(r,n)$. 
\end{thm}

\begin{proof}
The proof follows that of Proposition 3.1 in \cite{yan1994global}. Let $B := A + aI_{n}$ and consider the function $z(t) := \| I_{r} - U(t)^{\top} U(t) \|_{F}^{2}$, where $\| \bullet \|_{F}$ is the Frobenius norm. Its time derivative satisfies
\begin{align*}
\frac{d}{dt} z(t) &= -2 \mathrm{Tr} \left[ (I_{r} - U(t)^{\top} U(t) ) \frac{d}{dt}( U(t)^{\top} U(t)) \right]
\\
&=
-\frac{2}{\varepsilon} \mathrm{Tr} \Big{[}
U(t) (I_{r} - U(t)^{\top} U(t) ) ^{2} U(t)^{\top} (B^{\top} + B)
\Big{]}
\\
&\leq
-\frac{4}{\varepsilon} \mathrm{Tr} \Big{[}
U(t)^{\top} U(t) (I_{r} - U(t)^{\top} U(t) ) ^{2} \Big{]}
\lambda _{n}(B_{\mathrm{sym}})
\\
&\leq
-\frac{4}{\varepsilon}
\lambda _{r}( U(t)^{\top} U(t)) \lambda _{n}(B_{\mathrm{sym}})
z(t)
.
\end{align*}
By construction, $\lambda_{n}(B_{\rm sym}) > 0$. From Lemma \ref{lemma:svd_oja_solution}, we know that the smallest singular value of $U(t)$ is bounded below by $\alpha > 0$, which implies $\lambda_{r}(U(t)^{\top}U(t)) \geq \alpha^2 > 0$.
Thus, $z(t)$ satisfies 
\begin{align*}
\frac{d}{dt}z(t) \leq -\frac{4}{\varepsilon} \alpha ^{2} \lambda_{n}(B_{\rm sym}) z(t)
\end{align*}
which proves exponential convergence of $z(t)$ to 0. This result is equivalent to the exponential convergence of $U(t)$ to $\mathrm{St}(r,n)$.
\end{proof}

The convergence rate is tunable by the parameters $a$ and $\varepsilon$. A simple choice to ensure $A_{\rm sym} + aI_{n} > 0$ is to select $a > \|A\|_F = \sqrt{\mathrm{Tr}[A^{\top}A]}$. However, an exceedingly large value of $a$ makes the system dynamics stiff, requiring a smaller time step for stable numerical integration.

\begin{figure}[!htbp]
    	\centering
    	\includegraphics[keepaspectratio,width=\linewidth]{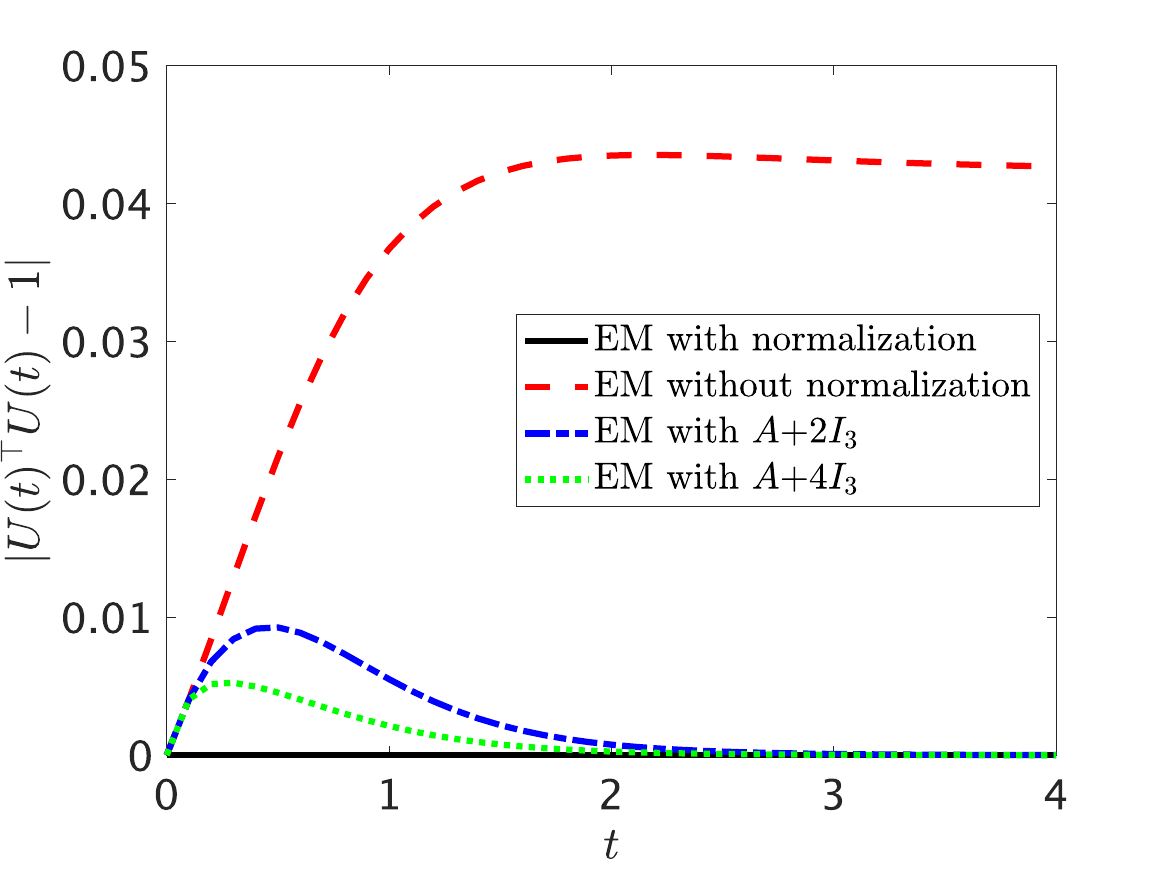}
    	\caption{Plot of $|U(t)^{\top} U(t) -1|$ for each time $t$ by the forward Euler method with normalization (black solid line), without normalization (red dashed line), with $A+2I_{3}$ (blue chain line), and with $A+4I_{3}$ (green dotted line) under the conditions in Example \ref{exmp:St13} with $U(0) = (\psi _{2} + \psi _{3})/\| \psi _{2} + \psi _{3} \|$. }
	\label{fig:example_adding_matrix}
	\end{figure}

	To demonstrate the stabilizing effect of the spectral shift $aI_n$ proposed in Theorem \ref{prop:positive_part_convergence}, we revisit the setting of Example \ref{exmp:St13}. The matrix $A$ in this example is not positive-definite. Consequently, simulating \eqref{eq:Oja_flow} with a standard forward Euler integrator (time step $h=0.1$) can cause the solution $U(t)$ to drift off the Stiefel manifold because of numerical errors, as illustrated by the red dashed line in Fig. \ref{fig:example_adding_matrix}.

In contrast, applying the same numerical scheme to the modified flow with $A+aI_n$ for $a=2$ and $a=4$ yields stable trajectories. 
Because $\lambda_3(A_{\mathrm{sym}}) \approx -1.47$, both values of $a$ satisfy the theorem's condition. As shown by the blue and green lines in Fig. \ref{fig:example_adding_matrix}, the solution may briefly leave the manifold because of initial numerical errors but is actively driven back. 
A larger value of $a$ results in faster convergence to the manifold; however, an excessively large $a$ can make the system stiff, requiring a prohibitively small time step $h$ for numerical stability. For comparison, the black solid line shows that periodic re-normalization also keeps the solution on the manifold. However, this approach incurs additional computational cost at each step; for instance, a Householder QR decomposition requires $O(nr^2)$ flops \cite[Algorithm 5.2.1]{golub2012matrix}, which can be demanding if $r$ is not small.

\begin{figure}[!htbp]
    	\centering
    	\includegraphics[keepaspectratio,width=\linewidth]{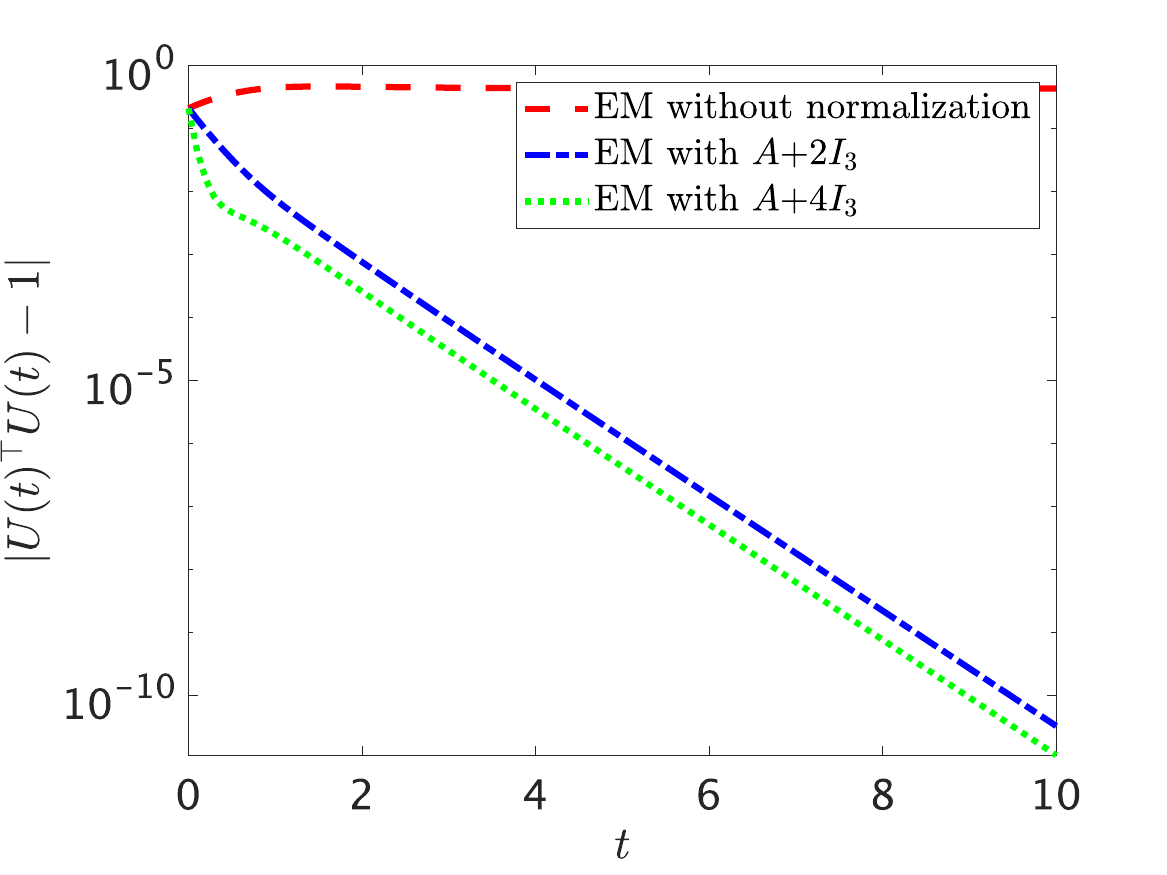}
    	\caption{Plot of $|U(t)^{\top} U(t) -1|$ for each time $t$ by the forward Euler method without normalization (red dashed line), with $A+2I_{3}$ (blue chain line), and with $A+4I_{3}$ (green dotted line) under the conditions in Example \ref{exmp:St13} with $U(0) = \frac{11}{10}(\psi _{2} + \psi _{3})/\| \psi _{2} + \psi _{3} \|$. }
	\label{fig:example_adding_matrix2}
	\end{figure}

Fig. \ref{fig:example_adding_matrix2} illustrates the case where the initial condition $U(0)$ starts outside the Stiefel manifold, plotted on a semi-logarithmic scale. The blue and green lines again demonstrate exponential convergence to the manifold. The trajectory with the larger value ($a=4$, green line) exhibits faster initial convergence. After $t=1$, the convergence rates appear similar, an artifact of the relatively large time step $h$. Using a smaller step size would show the superior convergence rate of the larger $a$ over a longer time interval.

	\begin{exmp}[Singular Subspace Extraction Algorithm]
Theorem \ref{prop:positive_part_convergence} can be used to improve the numerical stability of online singular subspace extraction algorithms. The Weingessel--Hornik (WH) algorithm \cite{weingessel1997svd} for extracting the singular subspaces of a matrix $A \in \mathbb{R}^{n \times m}$ is given by a set of coupled ordinary differential equations. This algorithm can be shown to be equivalent to the Oja flow for an augmented system \cite{fiori2003singular}:
\begin{align}
	\varepsilon \frac{d}{dt} X(t) = (I_{n+m} - X(t)X(t)^{\top}) \mathcal{A} X(t), \label{eq:SVD_Oja_flow}
\end{align}
where $X := [U^{\top}, V^{\top}]^{\top}/\sqrt{2}$ and $\mathcal{A} := \begin{bsmallmatrix} 0_{m} & A^{\top} \\ A & 0_{n} \end{bsmallmatrix}$.
Because $\mathcal{A}$ is symmetric but generally not positive-definite, a naive numerical implementation may be unstable if $X(t)$ deviates from the Stiefel manifold. 
To ensure numerical stability, we can apply Theorem \ref{prop:positive_part_convergence} by using the shifted matrix $\mathcal{B} := \mathcal{A} + aI_{n+m}$ for a sufficiently large $a>0$. 
This approach yields a modified WH algorithm with an additional term in each equation, which actively drives the solution back to the Stiefel manifold, eliminating the need for periodic re-normalization steps that can be computationally expensive. A detailed comparison with other truncated singular value decomposition algorithms is left for future work.
\end{exmp}

\subsection{Convergence of Oja Flow on the Stiefel Manifold}
\label{sec:convergence_on_SM}
		
	We then analyzed the convergence of the Oja flow on the Stiefel manifold.  
	The convergence problem can be reformulated by examining the dynamics of the projection matrix $P(t) := U(t)U(t)^{\top}$. If $U(t)$ solves the Oja flow \eqref{eq:Oja_flow}, then $P(t)$ evolves according to the matrix Riccati equation:
\begin{align}
 	\varepsilon \frac{d}{dt} P(t) = A P(t) + P(t) A^{\top} - P(t) (A+A^{\top}) P(t).
 	\label{eq:OjaRiccati}
\end{align}
To analyze this equation, we introduce an auxiliary linear system whose solution trajectory spans the same subspace:
\begin{align}
 	\varepsilon \frac{d}{dt} Z(t) = A Z(t), \quad Z(0) \in \mathbb{R}^{n\times r},
 	\label{eq:LinearDifferential_Eq}
\end{align}
where $Z(0)$ is assumed to be of full rank. The solution is $Z(t) = \mathrm{e}^{At/\varepsilon}Z(0)$, which preserves the initial rank. The projection matrix onto the subspace spanned by $Z(t)$ is $P'(t) := Z(t)(Z(t)^\top Z(t))^{-1}Z(t)^\top$. A direct calculation shows that $P'(t)$ also satisfies the Riccati equation \eqref{eq:OjaRiccati}. By setting $Z(0) = U(0)$, we ensure $P'(t) = P(t)$ for all $t$. This equivalence is powerful because $P'(t)$ has an explicit analytical form:
\begin{align*}
 	P'(t) =& \mathrm{e}^{At/\varepsilon} Z(0) \left( Z(0)^{\top} \mathrm{e}^{A^{\top}t/\varepsilon} \mathrm{e}^{At/\varepsilon} Z(0) \right)^{-1} 
	\\ &  \times
	Z(0)^{\top} \mathrm{e}^{A^{\top}t/\varepsilon}.
\end{align*}
This expression enables us to analyze the Oja flow's asymptotic behavior by studying the underlying linear system.

\begin{remark}
The trajectory $U(t)$ of the Oja flow is not, in general, a simple normalization of $Z(t)$. An additional time-varying rotation is involved, such that $U(t) = Z(t)(Z(t)^\top Z(t))^{-1/2} R(t)$ for some orthogonal matrix $R(t) \in \mathbb{R}^{r \times r}$. Characterizing this rotation is non-trivial, which motivates our analysis of the projection matrix $P(t)$ instead of $U(t)$ directly.
\end{remark}

\begin{thm}
\label{thm:convergence_Oja_flow}
Assume the eigenvalues of $A \in \mathbb{R}^{n\times n}$ satisfy the gap condition $\mathrm{Re}(\lambda_{m}(A)) > \mathrm{Re}(\lambda_{m+1}(A))$ for some integer $m$ with $r \leq m < n$.
Then, for any initial condition $U_0$ in the set
\begin{align*}
 	\mathcal{V}_{m,r} := \Bigg{\{} & \Psi(A) \begin{bmatrix} K_{m,r} \\ K_{m,\perp} \end{bmatrix} \in \mathrm{St}(r,n) \ \Bigg{|} \ K_{m,r} \in \mathbb{C}^{m\times r}, 
	\\ & \hspace{1cm}
	K_{m,\perp} \in \mathbb{C}^{(n-m)\times r}, \ \mathrm{rank}(K_{m,r}) = r
	 \Bigg{\}},
\end{align*}
the solution $U(t)$ of \eqref{eq:Oja_flow} converges exponentially to the invariant set
\begin{align*}
 	\mathcal{U}_{m,r} := \left\{ \Psi(A) \begin{bmatrix} K_{m,r} \\ 0_{n-m,r} \end{bmatrix} \in \mathrm{St}(r,n) \ \Bigg{|} \ K_{m,r} \in \mathbb{C}^{m\times r} \right\}.
\end{align*}
The convergence rate is governed by the exponent $-(\mathrm{Re}(\lambda_{m}) - \mathrm{Re}(\lambda_{m+1}) - \delta )/\varepsilon$, where $\delta >0$ is an arbitrarily small constant.
In the particular case where $r=m$, any solution starting in $\mathcal{V}_{r} := \mathcal{V}_{r,r}$ converges exponentially to an element of $\mathcal{U}_{r} := \mathcal{U}_{r,r}$.
\end{thm}

\begin{proof}
Let $\Lambda = \Psi^{-1}A\Psi$ be the Jordan normal form of $A$, where $\Psi = \Psi (A)$. We analyze the convergence of the projection matrix $P(t) = P'(t)= Z(t)(Z(t)^{\top}Z(t))^{-1}Z(t)^{\top}$, using its explicit form derived from the auxiliary system \eqref{eq:LinearDifferential_Eq} with $Z(0)=U_0 \in \mathcal{V}_{m,r}$. 

The key idea is to analyze the system in a scaled coordinate frame that factors out the dominant exponential growth. Let $Z'(t) = \mathrm{e}^{-\mathrm{Re}(\lambda_m)t/\varepsilon}Z(t)$. The projection matrix can be written as $P(t) = Z'(t)(Z'(t)^\top Z'(t))^{-1}Z'(t)^\top$. 
Then, $Z'(t)$ becomes 
\begin{align*}
Z^{\prime}(t) =& 
\Psi  \begin{bmatrix}
\mathrm{e}^{( \Lambda _{m} - \mathrm{Re} (\lambda _{m} ) I_{m} )t/\varepsilon}K_{m,r} 
\\ 
\mathrm{e}^{( \Lambda _{\perp} - \mathrm{Re} (\lambda _{m} ) I_{n-m})t/\varepsilon}K_{m,\perp}
\end{bmatrix}
\end{align*}
and $\Lambda_m \in \mathbb{C}^{m\times m}$ and $\Lambda_\perp \in \mathbb{C}^{(n-m) \times (n-m)}$ are Jordan blocks forming $\Lambda = \mbox{blk-diag}(\Lambda_m, \Lambda_\perp)$, corresponding to the dominant and non-dominant eigenvalues. 
Because of the eigenvalue gap condition, the components of the solution associated with $\Lambda_\perp$ decay exponentially relative to the components associated with $\Lambda_m$. Specifically, one can show that
\begin{align*}
    & \| \mathrm{e}^{\Lambda_{\perp} t/\varepsilon} K_{m,\perp} \| _{\rm ind}  \mathrm{e}^{-\mathrm{Re}(\lambda_m)t/\varepsilon} 
    \\
    \le &
    c \ \mathrm{e}^{(\mathrm{Re}(\lambda_{m+1}) - \mathrm{Re}(\lambda_m) + \delta)t/\varepsilon} \to 0,
\end{align*}
for an arbitrarily small $\delta > 0$ and a constant $c > 0$. 
Thus, $P(t) = P^{\prime}(t)$ converges to the subset $\{ \bar{U} \bar{U}^{\top} \ | \ \bar{U} \in \mathcal{U}_{m,r} \}$. 
This result proves that $U(t)$ converges to $\mathcal{U}_{m,r} $.  

For the special case $r=m$, the analysis can be refined to show exponential convergence to a specific projector, which implies convergence of $U(t)$ to the set $\mathcal{U}_r$. 
Because $K_{r} := K_{r,r}$ is square and full-rank, we can consider $P^{\prime}(t) = \hat{Z}(t) (\hat{Z}(t)^{\dagger} \hat{Z}(t) )^{-1}\hat{Z}(t)^{\dagger}$, where 
\begin{align*}
\hat{Z}(t) = Z(t) (\mathrm{e}^{ \Lambda _{r} t/\varepsilon}K_{r} )^{-1} 
=
\Psi  \begin{bmatrix} I_{r}
\\ 
\mathrm{e}^{\Lambda _{\perp} t/\varepsilon}K_{r,\perp}
K_{r}^{-1} \mathrm{e}^{ - \Lambda _{r} t/\varepsilon}
\end{bmatrix}
.
\end{align*}
Using norm inequality, 
\begin{align*}
& \| \mathrm{e}^{\Lambda _{\perp} t/\varepsilon}K_{r,\perp}
K_{r}^{-1} \mathrm{e}^{ - \Lambda _{r} t/\varepsilon} \| _{\rm ind}
\\
\leq &c^{\prime} \ \mathrm{e}^{(\mathrm{Re}(\lambda_{m+1}) - \mathrm{Re}(\lambda_m) + \delta)t/\varepsilon} 
\end{align*}
holds, where $\delta >0$ is an arbitrarily small constant and $c^{\prime}>0$.  
Thus, $\lim _{t \to \infty} \hat{Z}(t) = \Psi [I_{r}, 0_{r,n-r}]^{\top} $, proving that $U(t)$ converges to $\mathcal{U}_{r} $. 
Because any element of $\mathcal{U}_{r} $ is an equilibrium solution of \eqref{eq:Oja_flow}, $U(t)$ converges to an element of $\mathcal{U}_{r} $.  
\end{proof}

Theorem \ref{thm:convergence_Oja_flow} establishes convergence for the general case where $r \leq m$, thereby extending the recent result of Kuo et al. \cite{kuo2025asymptotic}, which addresses the specific case where $r=m$. 
A consequence of the proof is that the projection matrix $\bar{U}\bar{U}^\top$ is unique for any $\bar{U} \in \mathcal{U}_r$. The initial condition requires that the subspace spanned by $U_0$ have a non-trivial projection onto the dominant invariant subspace of $A$. As discussed in Section \ref{sec:basin}, this condition is not overly restrictive.

\begin{figure}[!htbp]
    	\centering
    	\includegraphics[keepaspectratio,width=8cm]{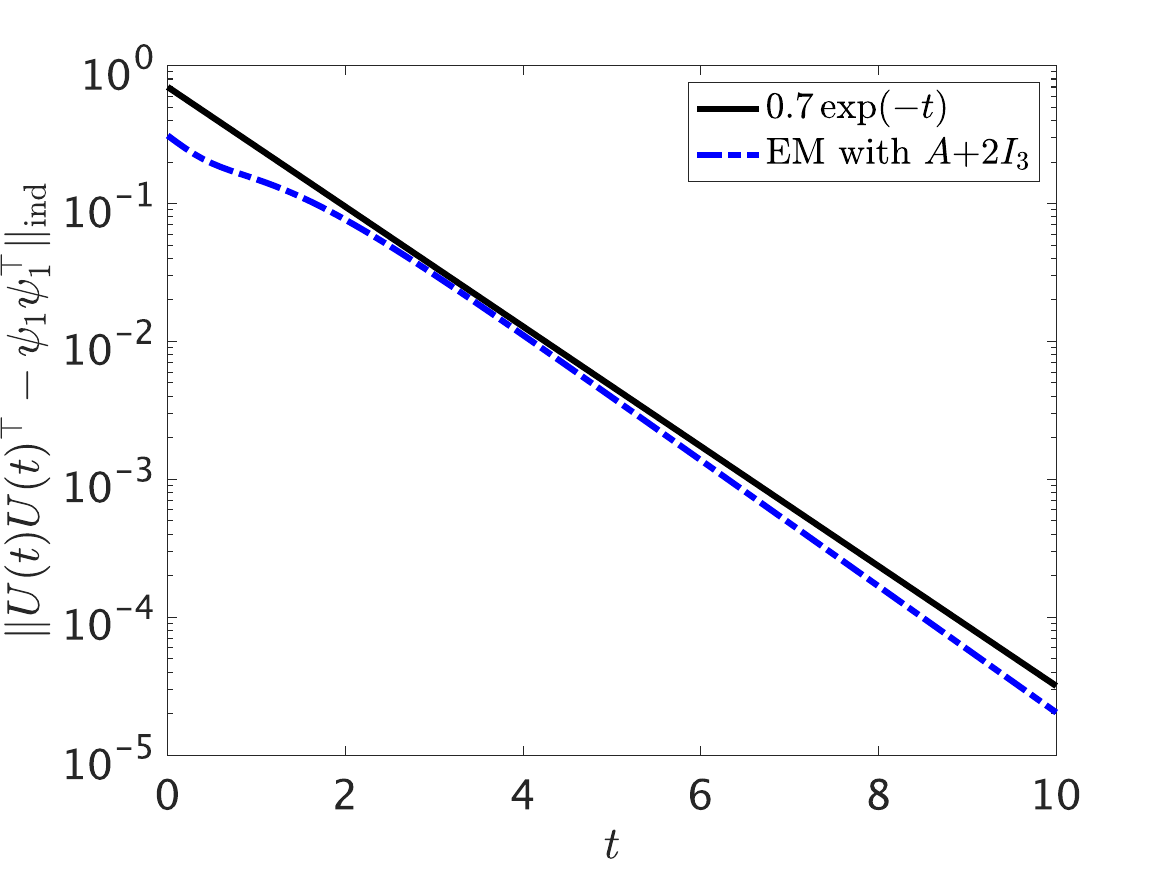}
    	\caption{Plot of $\| U(t) U(t)^{\top} - \psi _{1} \psi _{1}^{\top} \| _{\rm ind}$ for each time $t$ using the forward Euler method with $A+2I_{3}$ (blue chain line) under the conditions in Example \ref{exmp:St13} with $U(0) = (\psi _{1} + \psi _{2} + \psi _{3})/\| \psi _{1} + \psi _{2} + \psi _{3} \|$. The black line shows the upper bound of the convergence rate. }
	\label{fig_example_convergence_dominant_component_semilog}
	\end{figure}

Fig. \ref{fig_example_convergence_dominant_component_semilog} shows the convergence of the Oja flow \eqref{eq:Oja_flow} to the dominant mode $\mathcal{U}_{1}$. The matrix $A$ is defined in Example \ref{exmp:St13}, and the forward Euler scheme is used with time step $h=0.1$ and $\varepsilon =1$. The black line shows the upper bound of the convergence rate ($- (\mathrm{Re}(\lambda _{1}(A)) - \mathrm{Re}(\lambda _{2}(A))) = -1$). The figure shows that the convergence rate of $U(t)$ is slightly faster than the upper bound. To emphasize that the actual convergence rate is faster than the upper bound, we use $0.7 \mathrm{e}^{-t}$ rather than using the theoretical upper bound $\mathrm{e}^{-t}$.

\begin{remark}[On the case $r < m$]
When the desired subspace dimension $r$ is less than the number of dominant modes $m$, the eigenvalues of the projected matrix $\bar{U}^{\top}A\bar{U}$ for an equilibrium point $\bar{U} \in \mathcal{U}_{m,r}$ are not necessarily a subset of the eigenvalues of $A$. 
Notably, $\bar{U}$ is written as $\bar{U} = \Psi [K_{m,r}^{\top} , 0 _{r,n-m}]^{\top}$ and the QR decomposition yields $K_{m,r} = Q_{m,r}K_{r}$, where $Q_{m,r} \in \mathbb{C}^{m \times r}$ has orthonormal columns and $K_{r} \in \mathbb{C}^{r\times r}$ is invertible.  
Let $Q :=  \begin{bmatrix} Q_{m,r}^{\top} ,& 0_{r,(n-m)} \end{bmatrix}^{\top} \in \mathbb{C}^{n\times r}$ and $Q_{\perp} \in \mathbb{C}^{n \times (n-r)}$ so that $[Q,Q_{\perp} ]$ is a unitary matrix.  
Following similar arguments to the proof of \cite[Prop. 3]{TsuzukiOhki2024}, 
\begin{align*}
	& \bar{U}^{\top} A \bar{U} 
	\\
	=& 
	K_{r}^{\dagger} Q^{\dagger} \Psi ^{\dagger} \Psi 
	\\ &
	\times 
	\begin{bmatrix} QK_{r},& Q_{\perp}  \end{bmatrix}
	\begin{bmatrix} K_{r}^{-1}Q^{\dagger} \\ Q_{\perp}^{\dagger}  \end{bmatrix}
	\begin{bmatrix} \Lambda _{m} Q_{m,r}K_{r} \\ 0_{(n-r),r} \end{bmatrix}
	\\
	=&
	\begin{bmatrix}I_{r} , & K_{r}^{\dagger} Q^{\dagger} \Psi ^{\dagger} \Psi Q_{\perp} \end{bmatrix}
	\begin{bmatrix} K_{r}^{-1} Q_{m,r}^{\dagger} \Lambda _{m} Q_{m,r}K_{r} 
	\\ 
	Q_{\perp} ^{\dagger} [ ( \Lambda _{m} Q_{m,r}K_{r})^{\top}, 0_{(n-r),r} ^{\top} ]^{\top}
	\end{bmatrix}
	\\
	=&K_{r}^{-1} Q_{m,r}^{\dagger} \Lambda _{m} Q_{m,r}K_{r} 
	\\ & +
	K_{r}^{\dagger} Q^{\dagger} \Psi ^{\dagger} \Psi Q_{\perp} Q_{\perp}^{\dagger}  [ ( \Lambda _{m} Q_{m,r}K_{r})^{\top}, 0_{(n-r),r} ^{\top} ]^{\top}
	\end{align*}
holds. 
Consequently, the eigenvalues of $\bar{U}^{\top}A\bar{U}$ will generally not coincide with any of the eigenvalues of $A$ for $r<m$.
\end{remark}

Conversely, for the case where $r > m$, if the solution is initialized within a larger invariant subspace, it will converge while preserving its dominant components. This idea is formalized in the following discussion.

\begin{prop} \label{prop:r>m}
Assume there exist integers $m$ and $m'$ such that $1 \leq m < r \leq m' < n$ and that the eigenvalue gap conditions $\mathrm{Re}(\lambda_{m}) > \mathrm{Re}(\lambda_{m+1})$ and $\mathrm{Re}(\lambda_{m'}) > \mathrm{Re}(\lambda_{m'+1})$ both hold.
If the initial matrix $U(0)$ belongs to the set
\begin{align*}
    \Bigg{\{} \Psi(A) \begin{bmatrix} K_{m,r} \\ K_{m,\perp} \end{bmatrix} \in \mathcal{V}_{m',r} \mid &
    K_{m,r} \in \mathbb{C}^{m \times r},
    \\ & \mathrm{rank}(K_{m,r}) = m \Bigg{\}},
\end{align*}
then the solution $U(t)$ converges to an element of the set
\begin{align*}
    &\Bigg{\{} \Psi(A) \begin{bmatrix} K'_{m,r} \\ K'_{m'-m,\perp} \\ 0_{n-m',r} \end{bmatrix} \in \mathcal{U}_{m',r} \Bigg{|}  K'_{m,r} \in \mathbb{C}^{m \times r}, 
    \\ & \hspace{2cm}
    \mathrm{rank}(K'_{m,r}) = m, \ K'_{m'-m,\perp} \in \mathbb{C}^{(m'-m) \times r}
    \Bigg{\}}.
\end{align*}
\end{prop}

In essence, this proposition states that if the initial condition is contained within the $m'$-dominant invariant subspace and already spans the $m$-dominant invariant subspace, then the flow will converge to the $m'$-dominant subspace while preserving the span of the $m$-dominant component.

\begin{proof}
From the proof of Theorem \ref{thm:convergence_Oja_flow}, the projection matrix $P(t) = U(t)U(t)^\top$ converges to the set of projectors onto the subspace $\mathcal{U}_{m',r}$.
The asymptotic behavior of $P(t)$ is determined by the projection matrix $S(t) := V(t)(V(t)^\top V(t))^{-1}V(t)^\top$, where
\begin{align*}
    V(t) = \mathrm{e}^{At/\varepsilon} \Psi \begin{bmatrix}
        K_{m,r} \\ K_{m'-m,\perp} \\ 0_{n-m',r}
    \end{bmatrix}
    = \Psi \begin{bmatrix}
        \mathrm{e}^{\Lambda_m t/\varepsilon} K_{m,r} \\
        \mathrm{e}^{\Lambda_{m,m'} t/\varepsilon} K_{m'-m,\perp} \\
        0_{n-m',r}
    \end{bmatrix}.
\end{align*}
Here, $\Lambda_m \in \mathbb{C}^{m\times m}$ and $\Lambda_{m,m'} \in \mathbb{C}^{(m'-m)\times(m'-m)}$ are Jordan blocks containing the eigenvalues $\{\lambda_i\}_{i=1}^m$ and $\{\lambda_i\}_{i=m+1}^{m'}$, respectively.
The initial condition ensures that the matrix $[K_{m,r}^\top, K_{m'-m,\perp}^\top]^\top$ has full rank $r$, which guarantees that $V(t)^\top V(t)$ is invertible for all $t$.
Furthermore, because $K_{m,r}$ has full rank $m$ and $\mathrm{e}^{\Lambda_m t/\varepsilon}$ is non-singular, the block $\mathrm{e}^{\Lambda_m t/\varepsilon} K_{m,r}$ retains rank $m$ for all $t \ge 0$.
This result confirms that the component of the solution spanning the $m$-dominant subspace does not lose rank, which proves the proposition.
\end{proof}

\subsection{Evaluation of the Domain of Attraction}
\label{sec:basin}
	Next, we investigated the size of the domain of attraction, $\mathcal{V}_{m,r}$, within the Stiefel manifold.
	We consider only $m<n$; otherwise, $\mathcal{V}_{n,r} =  \mathrm{St}(r,n)$. 
Intuitively, an element $U \in \mathrm{St}(r,n)$ is determined by its coordinate representation $K = \Psi^{-1}U$. The condition for convergence, as stated in Theorem \ref{thm:convergence_Oja_flow}, is that a specific sub-block of $K$, denoted $K_{m,r}$, must have full rank. For a randomly chosen $U$, it is highly improbable that this sub-block would be rank-deficient. This result suggests that the set of initial conditions leading to convergence, $\mathcal{V}_{m,r}$, should encompass almost the entire manifold. We formalize this claim in the following propositions.

\begin{prop}\label{prop:volume_basin}
The set $\mathcal{V}_{m,r}$ is dense in $\mathrm{St}(r,n)$. That is, its closure is the entire manifold: $\overline{\mathcal{V}_{m,r}} = \mathrm{St}(r,n)$.
\end{prop}

\begin{proof}
To prove that $\mathcal{V}_{m,r}$ is dense in $\mathrm{St}(r,n)$, we must show that for any $U \in \mathrm{St}(r,n)$ and any open neighborhood $\mathcal{N}_{\delta}(U)$ of $U$, the intersection $\mathcal{N}_{\delta}(U) \cap \mathcal{V}_{m,r}$ is non-empty \cite[Thm 17.5]{munkres2018topology}. We define a neighborhood as
\begin{align}
    \mathcal{N}_{\delta}(U) := \{ V \in \mathrm{St}(r,n) \mid \| U - V \|_{\rm ind} < \delta \}, \quad \delta > 0.
    \label{eq:neighborhood}
\end{align}
If $U \in \mathcal{V}_{m,r}$, the condition is trivially satisfied. Therefore, we consider the case where $U \in \mathrm{St}(r,n) \setminus \mathcal{V}_{m,r}$, which implies that its coordinate representation $K = \Psi^{-1}U = [K_{m,r}^\top, K_{m,\perp}^\top]^\top$ has a rank-deficient sub-block, i.e., $\mathrm{rank}(K_{m,r}) < r$.

Notably, for any orthogonal matrix $Q \in \mathbb{R}^{n\times n}$ and any $U \in \mathrm{St}(r,n)$, $Q U \in \mathrm{St}(r,n)$ holds. Consider a skew-symmetric matrix $W = -W^{\top} \in \mathbb{R}^{n\times n}$ with $\| W \| _{\mathrm{ind}} =1$ and $Q = \exp (\bar{\delta} W)$ for a small $\bar{\delta}>0$. We show that we can choose $\bar{\delta}>0$ such that $QU \in \mathcal{N}_{\delta} (U)$ for a given $\delta >0$. Using the operator norm inequality yields
\begin{align*}
\| U - QU \| _{\mathrm{ind}}
&\leq \| I_{n} - \exp (\bar{\delta} W) \| _{\mathrm{ind}} \| U \| _{\mathrm{ind}}
\\
&= \| I_{n} - \exp (\bar{\delta} W) \| _{\mathrm{ind}}
\end{align*}
and because $\| W \| _{\mathrm{ind}} =1$,
\begin{align*}
\| I_{n} - \exp (\bar{\delta} W) \| _{\mathrm{ind}}
&\leq
\sum _{k=1} ^{\infty} \frac{ \bar{\delta} ^{k} }{k!} = e^{ \bar{\delta}} -1
\end{align*}
holds. Taking $\bar{\delta}>0$ such that $\bar{\delta} < \ln \left( \delta +1 \right)$ ensures $QU \in \mathcal{N}_{\delta} (U)$.

Our strategy is to construct a matrix $V$ that is arbitrarily close to $U$ but lies in $\mathcal{V}_{m,r}$. We construct $V$ by applying a small orthogonal rotation to $U$, i.e., $V = QU$, where $Q = \exp(\bar{\delta}W)$ for a small scalar $\bar{\delta} > 0$ and a skew-symmetric matrix $W$ with $\|W\|_{\rm ind} = 1$. 

The crucial step is to choose $W$ such that the resulting matrix $V$ is in $\mathcal{V}_{m,r}$. Let $K'_{m,r} \in \mathbb{C}^{m\times r}$ be a matrix of full rank $r$ such that $\Psi_m K'_{m,r} \in \mathrm{St}(r,n)$, where $\Psi _m = [\psi_1 ,\dots , \psi _{m}] \in \mathbb{C}^{n\times m}$. 
We construct the rotation generator as $W = c(X Y^\dagger - Y X^\dagger)$, where $X = \Psi_m K'_{m,r} \in \mathrm{St}(r,n)$, $Y = U = \Psi K$, and $c$ is a normalization constant. 
Notably, $\mathrm{rank}(K'_{m,r}) =r$. 
This choice of $W$ perturbs $U$ in a direction that mixes its components with those of a full-rank element.

The coordinate representation of the perturbed matrix is $\Psi^{-1}V = \Psi^{-1}QU$. A first-order Taylor expansion yields
\begin{align*}
    \Psi^{-1}V =& \Psi^{-1}\exp(\bar{\delta}W)\Psi K = K + \bar{\delta}\Psi^{-1}W\Psi K + O(\bar{\delta}^2)
    \\
    =&
    K + c\bar{\delta}
    \left(
    \begin{bmatrix}
    K'_{m,r} \\ 0_{n-m,r}
    \end{bmatrix}
    -
    K (\Psi _{m} K'_{m,r})^{\dagger} U
    \right)
    + O(\bar{\delta}^2)
    .
\end{align*}
The top $m \times r$ block of the perturbed coordinates, which we denote as $K'_{m,r}(\bar{\delta})$, is
\begin{align*}
    K'_{m,r}(\bar{\delta}) = K_{m,r} (I_{r} -  c\bar{\delta} (X^\dagger U)  ) + c\bar{\delta} K'_{m,r} + O(\bar{\delta}^2).
\end{align*}
Because $K_{m,r}$ is rank-deficient and $K'_{m,r}$ is full-rank, the first-order term $c\bar{\delta}K'_{m,r}$ perturbs $K_{m,r}$ in a direction that restores its rank. For a sufficiently small $\bar{\delta} > 0$, the higher-order terms do not alter this outcome, so $\mathrm{rank}(K'_{m,r}(\bar{\delta})) = r$. This result implies that $V=QU \in \mathcal{V}_{m,r}$.
Because we can find such a $V$ in any arbitrarily small neighborhood of any $U$, we conclude that $\mathcal{V}_{m,r}$ is dense in $\mathrm{St}(r,n)$.
\end{proof}

Proposition \ref{prop:volume_basin} establishes that $\mathcal{V}_{m,r}$ is dense, but this result is not sufficient to make claims about its volume. For example, the set of rational numbers is dense in the real line but has zero Lebesgue measure. 
On the Stiefel manifold, there exists a uniform invariant measure \cite[\S 1.4.3]{Chikuse2003}. 
To show that $\mathcal{V}_{m,r}$ has the same volume as $\mathrm{St}(r,n)$, the following proposition fills the gap.

\begin{prop}\label{prop:connectness}
For any $U \in \mathcal{V}_{m,r}$,
\begin{enumerate}
\item $\mathcal{N}_{\delta}(U)$ is a path-connected set for $\delta \in (0,2)$.
\item There exists $\delta >0$ such that $\mathcal{N}_{\delta}(U) \subset \mathcal{V}_{m,r}$, i.e., $ \mathcal{V}_{m,r}$ is an open subset of $\mathrm{St}(r,n)$.
\end{enumerate}
\end{prop}

\begin{proof}
1) Path-connectedness of neighborhoods: 
For any $U \in \mathcal{V}_{m,r}$, the singular value decomposition is $U = Q [I_{r} , 0_{r,n-r}]^{\top}$, where $Q \in \mathbb{R}^{n\times n}$ is an orthogonal matrix. 
Hence, for any $ U_{1} , U_{2} \in \mathrm{St}(r,n)$, there exist orthogonal matrices $Q\in \mathbb{R}^{n\times n}$ such that $U_{2} = Q U_{1}$. 
	Hence, $V \in \mathcal{N}_{\delta}(U)$ is represented as $V=QU$, and from the definition of $\mathcal{N}_{\delta} (U)$, we can take $Q$ that satisfies 
	\begin{align*}
	\| U - QU \| _{\rm ind} 
	\leq & \| I_{n} - Q \| _{\rm ind} 
	= \max _{i=1,\dots ,n} | 1 - \lambda _{i}(Q) |  < \delta .
	\end{align*}
	Because $Q \in \mathbb{R}^{n\times n}$ is orthogonal, all the eigenvalues are on the unit circle within the $\delta$-ball from $1$, and therefore, $Q$ does not contain its eigenvalue on $-1$ if we take $\delta <2$. 
	This result means that $\mathrm{det}(Q) =1$, i.e., $Q$ is in a special orthogonal group $SO(n)$.  
	Because $SO(n)$ is a simple Lie group, $\mathcal{N}_{\delta} (U)$ is path-connected for small $\delta >0$.

	2) Openness of $\mathcal{V}_{m,r}$: We must show that for any $U \in \mathcal{V}_{m,r}$, there exists a neighborhood $\mathcal{N}_{\delta}(U)$ that is entirely contained within $\mathcal{V}_{m,r}$.
	If $U \in \mathcal{V}_{m,r}$, its coordinate block $K_{m,r}$ has full rank. 
	The rank of a full-rank matrix does not decrease its rank under small perturbations. 
	Therefore, for any sufficiently small perturbation of $U$ to a nearby point $V \in \mathcal{N}_{\delta}(U)$, the corresponding coordinate block $K'_{m,r}$ will also have full rank. This result implies there exists a $\delta > 0$ such that $\mathcal{N}_{\delta}(U) \subset \mathcal{V}_{m,r}$, proving that $\mathcal{V}_{m,r}$ is an open set.
\end{proof}

	The Stiefel manifold $\mathrm{St}(r,n)$ is a compact manifold equipped with a uniform invariant measure, giving it a finite, well-defined volume \cite[\S 1.4.4]{Chikuse2003}. 
	Because $\mathcal{V}_{m,r}$ is an open and dense subset of $\mathrm{St}(r,n)$, its complement must be a closed set with an empty interior, which has measure zero. This leads to our main result on the domain of attraction. 

\begin{thm}\label{thm:volume_of_basin}
The volume of the set $\mathcal{V}_{m,r}$ is equal to the volume of the Stiefel manifold $\mathrm{St}(r,n)$.
\end{thm}

A consequence of Theorems \ref{thm:convergence_Oja_flow} and \ref{thm:volume_of_basin} implies that there is no other stable invariant set for the Oja flow \eqref{eq:Oja_flow}. 
Theorem \ref{thm:volume_of_basin} implies that if an initial condition $U(0)$ is chosen uniformly at random from $\mathrm{St}(r,n)$, it will belong to the domain of attraction $\mathcal{V}_{m,r}$ with probability one. Combining this result with those described in Section \ref{sec:Oja_flow_Euclidean}, one can confidently initialize the Oja flow from an arbitrary full-rank matrix in $\mathbb{R}^{n\times r}$. By using a suitable shift $a>0$ to ensure $A_{\rm sym}+aI_n > 0$, the trajectory will first converge to the Stiefel manifold and then, with probability one, proceed to converge to the dominant invariant subspace $\mathcal{U}_{m,r}$.

\subsection{Change in the Number of Dominant Component Subspaces}
	
	In this analysis, we assume that a solution $\bar{U}_r \in \mathcal{U}_r$ spanning an $r$-dimensional dominant invariant subspace has been found. We now explore efficient methods for obtaining related subspaces of either higher or lower dimension.

\subsubsection{Expanding the Subspace}

The following proposition provides a computationally efficient method for finding the principal components corresponding to a larger, $(r+\ell)$-dimensional subspace, given the solution for the $r$-dimensional one.

\begin{prop}[Subspace Expansion]
\label{prop:subspace_expansion}
Assume the eigenvalues of $A \in \mathbb{R}^{n\times n}$ satisfy the gap conditions $\mathrm{Re}(\lambda_{r}(A)) > \mathrm{Re}(\lambda_{r+1}(A))$ and $\mathrm{Re}(\lambda_{r+\ell}(A)) > \mathrm{Re}(\lambda_{r+\ell+1}(A))$ for $r \geq 1$ and $1 \leq \ell < n-r$.
Given a solution $\bar{U}_r \in \mathcal{U}_r$, let $\bar{U}_\perp \in \mathrm{St}(n-r, n)$ be an orthonormal basis for the orthogonal complement of the subspace spanned by $\bar{U}_r$, such that $\bar{U}_\perp \bar{U}_\perp^\top = I_n - \bar{U}_r \bar{U}_r^\top$. Consider the reduced-order Oja flow:
\begin{align}
    \varepsilon \frac{d}{dt} u(t) = (I_{n-r} - u(t)u(t)^\top) A_{\bar{U}_\perp} u(t),
    \label{eq:reduced_Oja_flow}
\end{align}
where $u(t) \in \mathbb{R}^{(n-r) \times \ell}$ and $A_{\bar{U}_\perp} := \bar{U}_{\perp}^{\top} A \bar{U}_{\perp}$.
Then, for an initial condition $u(0)$ such that $[\bar{U}_r, \bar{U}_\perp u(0)] \in \mathcal{V}_{r+\ell}$, the combined solution $U_{r+\ell}(t) := [\bar{U}_r, \bar{U}_\perp u(t)]$ converges to an element of $\mathcal{U}_{r+\ell}$.
\end{prop}

\begin{proof}
Because $\bar{U}_r$ spans an invariant subspace of $A$, the matrix $A$ is block upper-triangular in the basis $[\bar{U}_r, \bar{U}_\perp]$, meaning $\bar{U}_\perp^\top A \bar{U}_r = 0$. Consequently, the eigenvalues of the projected matrix $A_{\bar{U}_\perp}$ are precisely the remaining $n-r$ eigenvalues of $A$, i.e., $\{\lambda_{r+1}(A), \dots, \lambda_n(A)\}$.
From Theorem \ref{thm:convergence_Oja_flow}, the solution $u(t)$ of the reduced-order Oja flow \eqref{eq:reduced_Oja_flow} converges to the $\ell$-dimensional dominant invariant subspace of $A_{\bar{U}_\perp}$.

Here, we verify that the dynamics of the composite matrix $U_{r+\ell}(t) := [\bar{U}_r, \bar{U}_\perp u(t)]$ are equivalent to the full $(r+\ell)$-dimensional Oja flow. The time derivative is $\frac{d}{dt}U_{r+\ell}(t) = [0_{n,r}, \bar{U}_\perp \frac{d}{dt}u(t)]$. Substituting the dynamics from \eqref{eq:reduced_Oja_flow} yields
\begin{align*}
    \frac{d}{dt}U_{r+\ell}(t) = \frac{1}{\varepsilon} \left[ 0_{n,r}, \bar{U}_\perp (I_{n-r} - u(t)u(t)^\top) A_{\bar{U}_\perp} u(t) \right].
\end{align*}
Conversely, the full Oja flow for $U_{r+\ell}(t)$ is $(I_n - U_{r+\ell}(t)U_{r+\ell}(t)^\top) A U_{r+\ell}(t)$. Using the facts that $\bar{U}_\perp^\top A \bar{U}_r = 0$ and $I_n - U_{r+\ell}(t) U_{r+\ell}(t)^\top = \bar{U}_\perp(I_{n-r} - u(t)u(t)^\top)\bar{U}_\perp^\top$, one can show that this expression simplifies to the same result. Therefore, the trajectory of $U_{r+\ell}(t)$ is identical to that of an $(r+\ell)$-dimensional Oja flow, and its convergence to $\mathcal{U}_{r+\ell}$ is guaranteed by Theorem \ref{thm:convergence_Oja_flow}.
\end{proof}

Because \eqref{eq:reduced_Oja_flow} is a lower-dimensional Oja flow, the stabilization technique from Theorem \ref{prop:positive_part_convergence} can be applied to ensure robust numerical computation. Although constructing the basis $\bar{U}_\perp$ for the orthogonal complement may be non-trivial, it can often be performed as an offline computation.

\subsubsection{Reducing the Subspace}

Next, we consider methods for extracting a lower-dimensional, $\tilde{r}$-dominant subspace (with $\tilde{r} < r$) from a given solution $\bar{U}_r \in \mathcal{U}_r$.

\paragraph{Method 1: Eigendecomposition}
As shown in \cite[Prop. 3]{TsuzukiOhki2024}, the projected matrix $\bar{U}_r^\top A \bar{U}_r$ is related to the dominant Jordan block of $A$ by a similarity transform: $\bar{U}_r^\top A \bar{U}_r = K_r^{-1}\Lambda_r K_r$. The matrix of eigenvectors of this small $r \times r$ problem gives the coordinate transformation $K_r$. The first $r$ eigenvectors of the original matrix $A$ can then be recovered as $\Psi_r = \bar{U}_r K_r$.
To obtain the $\tilde{r}$-dominant subspace, one can simply select the first $\tilde{r}$ columns of $\Psi_r$ to form $\Psi_{\tilde{r}}$, and then re-orthonormalize to obtain $\bar{U}_{\tilde{r}} \in \mathcal{U}_{\tilde{r}}$. This approach is highly efficient for small $r$, but the cost of the $r \times r$ eigendecomposition can become significant for larger $r$.

\paragraph{Method 2: Recursive Oja Flow}
An alternative iterative approach is to apply the Oja flow recursively. Because the small matrix $A_{\bar{U}_r} := \bar{U}_r^\top A \bar{U}_r \in \mathbb{R}^{r \times r}$ is readily available, we can solve the following reduced-order Oja flow on its domain:
\begin{align*}
    \varepsilon \frac{d}{dt} \tilde{U}(t) = (I_r - \tilde{U}(t)\tilde{U}(t)^\top) A_{\bar{U}_r} \tilde{U}(t),
\end{align*}
where $\tilde{U}(0) \in \mathrm{St}(\tilde{r}, r)$ is chosen from the domain of attraction for this smaller system. The solution $\tilde{U}(t)$ converges to an element $\tilde{U}_\infty \in \mathcal{U}_{\tilde{r}}^r$, where $\mathcal{U}_{\tilde{r}}^r$ is the set of $\tilde{r}$-dominant invariant subspaces of $A_{\bar{U}_r}$.

The key insight is that the solution of this small-scale problem directly provides the projection needed to extract the desired subspace from the original solution $\bar{U}_r$. An element of $\mathcal{U}_{\tilde{r}}^r$ has the form $K_r^{-1} [\tilde{K}_{\tilde{r}}^\top, 0]^\top$. Therefore, the desired $\tilde{r}$-dimensional subspace of the original system is obtained by the simple matrix product:
\begin{align*}
    \bar{U}_{\tilde{r}} := \bar{U}_r \tilde{U}_\infty \in \mathcal{U}_{\tilde{r}}.
\end{align*}
This method avoids a full eigendecomposition and can be advantageous when $r$ is moderately large.

	\section{Applications to Control Problems}
	\label{sec:applications}

	The extraction of dominant eigenvalues and eigenvectors is essential for analyzing large-scale dynamical systems and recursive algorithms. We demonstrate how to utilize the Oja flow \eqref{eq:Oja_flow} for model reduction, low-rank controller synthesis, and singularly perturbed systems. We also provide some theoretical guarantees for these applications.

	\subsection{Model Reduction for Linear Dynamical Systems}

	In practical control problems, model reduction is crucial for designing controllers for large and complex systems \cite{obinata2000model,antoulas2005overview,benner2015survey,SnowdenGraafTindall2017}.
Many existing methods, such as balanced truncation, primarily focus on stable, linear, time-invariant systems \cite{antoulas2005approximation,gugercin2008h_2,sato2018structure}.
Balanced truncation relies on computing controllability and observability Gramians, a task that can be computationally prohibitive for large-scale systems. One of its key advantages is that the modeling error can be bounded in the $H^\infty$ or $H^2$ norm.
Although extensions such as time-limited and linear--quadratic--Gaussian (LQG) balanced truncation exist for unstable systems \cite{zhou1999balanced,opdenacker1985lqg,goyal2019time,das2022h}, they still require solving large-scale algebraic Riccati equations to obtain the necessary Gramians.
Other prominent techniques, such as Hankel norm approximation and Krylov subspace methods, are also often computationally demanding or are restricted to stable systems \cite{kung1981optimal,antoulas2005approximation}.
In this section, we propose model reduction methods and a related controller synthesis framework based on the Oja flow, which avoids the need to solve for Gramians directly.

For convenience, given a matrix $Y \in \mathbb{R}^{n \times q}$, we adopt the standard notation for projected systems: $A_Y := Y^\top A Y$, $B_Y := Y^\top B$, and $C_Y := CY$.

\subsubsection{Properties of Oja Flow Associated with $A$ and $A^{\top}$}

Assume that for a given matrix $A \in \mathbb{R}^{n \times n}$, the eigenvalue gap condition $\mathrm{Re}(\lambda_r(A)) > \mathrm{Re}(\lambda_{r+1}(A))$ holds.
The Oja flow \eqref{eq:Oja_flow} applied to $A$ and $A^\top$ can then extract the dominant right and left invariant subspaces, respectively. We denote the resulting steady-state solutions as $\bar{U} \in \mathcal{U}_r(A)$ and $\bar{V} \in \mathcal{U}_r(A^\top)$, where $\mathcal{U}_r(X)$ is the stable equilibrium set for the Oja flow associated with matrix $X$.
Let $\bar{U}_\perp \in \mathrm{St}(n-r, n)$ and $\bar{V}_\perp \in \mathrm{St}(n-r, n)$ be orthonormal bases for the orthogonal complements of the subspaces spanned by $\bar{U}$ and $\bar{V}$, respectively. Using the orthogonal matrices $Q_{\bar{U}} := [\bar{U}, \bar{U}_\perp]$ and $Q_{\bar{V}} := [\bar{V}, \bar{V}_\perp]$, we can subject $A$ to similarity transformations that reveal a block-triangular structure:
\begin{align}
	Q_{\bar{U}}^{\top} A Q_{\bar{U}} &= 
	\begin{bmatrix}
	A_{\bar{U}} & \bar{U}^{\top} A \bar{U}_{\perp}
	\\
	0_{n-r,r} & A_{\bar{U}_\perp}
	\end{bmatrix},
	\label{eq:blk_upper_triangle}
	\\
	Q_{\bar{V}}^{\top} A Q_{\bar{V}} &=
	\begin{bmatrix}
	A_{\bar{V}} & 0_{r,n-r} 
	\\
	\bar{V}_{\perp}^{\top} A\bar{V}  & A_{\bar{V}_\perp}
	\end{bmatrix}.
	\label{eq:blk_lower_triangle}
\end{align}
The $(2,1)$ block of \eqref{eq:blk_upper_triangle} is zero because $\bar{U}$ is an equilibrium point of its Oja flow, which implies $(I_n - \bar{U}\bar{U}^\top)A\bar{U} = 0_{n,r}$, and therefore $\bar{U}_\perp^\top A \bar{U} = 0_{n-r,r}$. Similarly, the $(1,2)$ block of \eqref{eq:blk_lower_triangle} is zero, which means $\bar{V}^{\top} A \bar{V}_{\perp} = 0_{r,n-r}$. 

The quadruple $(\bar{U}, \bar{U}_\perp, \bar{V}, \bar{V}_\perp)$ exhibits the following important properties.

\begin{lemma} \label{lem:properties_U_Uperp}
The matrix $\bar{U}_\perp$ spans the minor left invariant subspace of $A$; specifically, $\bar{U}_{\perp} \in \mathcal{U}_{n-r}(-A^{\top})$. Similarly, $\bar{V}_{\perp} \in \mathcal{U}_{n-r}(-A)$. Furthermore, the matrix $\bar{V}^{\top} \bar{U}$ is non-singular.
\end{lemma}
\begin{proof}
Any $\bar{U} \in \mathcal{U}_{r}(A)$ has the representation
\begin{align*}
\bar{U} = \Psi (A) \begin{bmatrix}
K_{r} \\ 0_{n-r,r}
\end{bmatrix}
,\ K_{r} \in \mathbb{C}^{r\times r},\
\mathrm{det}(K_{r}) \neq 0
.
\end{align*}
From the definition, $\bar{U}_{\perp}^{\top} \bar{U} = 0_{n-r,r}$. The columns of $\bar{U}_{\perp}$ span the complementary subspace of the dominant $r$-dimensional eigensubspace. Thus, $\bar{U}_{\perp}$ is associated with the remaining $n-r$ eigenvalues, and its form is
\begin{align*}
\bar{U}_{\perp} = \Psi (A)^{-\dagger} \begin{bmatrix}
0_{r,n-r} \\ \tilde{K}_{\perp}
\end{bmatrix}
,
&\ \tilde{K}_{\perp} \in \mathbb{C}^{(n-r) \times (n-r)},\\
& \mathrm{det}(\tilde{K}_{\perp}) \neq 0
.
\end{align*}
Because $A = \Psi (A) \Lambda \Psi (A)^{-1}$, we have
\begin{align*}
A^{\top} =& \Psi (A)^{-\dagger} \Lambda ^{\dagger} \Psi (A)^{\dagger}
= \Psi (A^{\top}) \Lambda \Psi (A^{\top}) ^{-1},
\end{align*}
where $\Psi (A^{\top}) = \Psi (A)^{-\dagger} P$ and $P \in \mathbb{R}^{n\times n}$ is a permutation matrix.
From the eigenvalue gap condition, $P = \mbox{blk-diag}(P_{r}, P_{\perp})$, where $P_{r} \in\mathbb{R}^{r\times r}$ and $P_{\perp} \in \mathbb{R}^{(n-r) \times (n-r)}$ are permutation matrices corresponding to the dominant eigenvalues; $\bar{V} \in \mathcal{U}_{r}(A^{\top})$ is given by
\begin{align*}
\bar{V} = \Psi (A^{\top}) \begin{bmatrix}
L_{r}^{\prime} \\ 0_{n-r,r}
\end{bmatrix}
=
\Psi (A)^{-\dagger} 	\begin{bmatrix}
L_{r} \\ 0_{n-r,r}
\end{bmatrix},
\end{align*}
where $L_{r}^{\prime} \in \mathbb{C}^{r\times r}$ is non-singular and $L_{r} = P_{r} L_{r}^{\prime}$. 
From the definition of $\mathcal{U}_{r}(\bullet )$, $\bar{U}_{\perp} \in \mathcal{U}_{n-r}(-A^{\top})$ and $\bar{V}_{\perp} \in \mathcal{U}_{n-r}(-A)$. 

Using these representations, the product $\bar{V}^{\top} \bar{U}$ is
\begin{align*}
\bar{V}^{\top} \bar{U} &= \left( \Psi (A)^{-\dagger} \begin{bmatrix}
L_{r} \\ 0_{n-r,r}
\end{bmatrix} \right)^{\top} \left( \Psi (A) \begin{bmatrix}
K_{r} \\ 0_{n-r,r}
\end{bmatrix} \right)
\\
&= \begin{bmatrix}
L_{r}^{\dagger} & 0
\end{bmatrix} \Psi (A)^{-1} \Psi (A) \begin{bmatrix}
K_{r} \\ 0
\end{bmatrix}
= L_{r}^{\dagger} K_{r}
.
\end{align*}
Because $\mathrm{det}(L_r) \neq 0$ and $\mathrm{det}(K_r) \neq 0$, the matrix $\bar{V}^{\top} \bar{U}$ is non-singular.
\end{proof}

According to Lemma \ref{lem:properties_U_Uperp}, the Oja flow can be expressed as a set of coupled equations for a basis and its orthogonal complement:
\begin{align*}
	\varepsilon \frac{d}{dt} U(t) &= U_{\perp}(t) U_{\perp}(t)^{\top} A U(t), \\
	\varepsilon \frac{d}{dt} U_{\perp}(t) &= -U(t) U(t)^{\top} A^{\top} U_{\perp}(t).
\end{align*}

We now relate the two reduced-order matrices, $A_{\bar{U}}$ and $A_{\bar{V}}$.

\begin{lemma}\label{lem:similar_transformation}
The matrices $A_{\bar{U}}$ and $A_{\bar{V}}$ are related by a similarity transformation:
\begin{align*}
	A_{\bar{U}} = (\bar{V}^{\top} \bar{U})^{-1} A_{\bar{V}} (\bar{V}^{\top} \bar{U}).
\end{align*}
Consequently, they share the same eigenvalues, which are the $r$ dominant eigenvalues of $A$.
\end{lemma}
\begin{proof}
From their coordinate representations, we know that $A_{\bar{U}} = K_r^{-1}\Lambda_r K_r$ and $A_{\bar{V}} = L_r^{\dagger}\Lambda_r L_r^{-\dagger}$. From the proof of Lemma \ref{lem:properties_U_Uperp}, we have $\bar{V}^\top \bar{U} = L_r^\dagger K_r$. A direct substitution shows that the similarity transformation holds.
\end{proof}

\subsubsection{Model Reduction and its Properties}

The Oja flow extracts the invariant subspace corresponding to the $r$ eigenvalues with the largest real parts. This property renders it a natural tool for identifying and isolating the dominant dynamics of a linear time-invariant (LTI) system, which is the cornerstone of model reduction. Consider the LTI system:
\begin{align}
	\frac{d}{dt}x(t) = Ax(t) + Bu(t), \quad y(t) = Cx(t),
	\label{eq:system}
\end{align}
where $x \in \mathbb{R}^n$, $u \in \mathbb{R}^m$, and $y \in \mathbb{R}^p$ are the state, input, and output, respectively. Throughout the analysis discussed in this section, we assume the eigenvalue gap condition $\mathrm{Re}(\lambda_r(A)) > \mathrm{Re}(\lambda_{r+1}(A))$ holds for a fixed dimension $r < n$.

Projecting the system \eqref{eq:system} onto the subspace spanned by $\bar{U} \in \mathcal{U}_r(A)$ yields a reduced-order model $(A_{\bar{U}}, B_{\bar{U}}, C_{\bar{U}})$. Although the basis $\bar{U}$ is unique only up to an orthogonal transformation, the resulting input--output behavior of the reduced model is unique. This result can be seen by examining the reduced model's impulse response, $C_{\bar{U}} \mathrm{e}^{A_{\bar{U}}t} B_{\bar{U}}$. Using the property that for an invariant subspace, $\mathrm{e}^{At}\bar{U} = \bar{U}\mathrm{e}^{A_{\bar{U}}t}$, this result becomes
\begin{align*}
	C_{\bar{U}} \mathrm{e}^{A_{\bar{U}}t} B_{\bar{U}} = C \bar{U} \mathrm{e}^{A_{\bar{U}}t} \bar{U}^\top B = C \mathrm{e}^{At} \bar{U}\bar{U}^\top B,
\end{align*}
which implies the systems $(A_{\bar{U}}, B_{\bar{U}}, C_{\bar{U}})$ and $(A, \bar{U}\bar{U}^{\top} B, C)$ are identical. 
Because the projector $\bar{U}\bar{U}^\top$ is unique for the subspace $\mathcal{U}_r(A)$, the impulse response is independent of the specific choice of basis $\bar{U}$. Because the eigenvalues of $A_{\bar{U}}$ are the $r$ dominant eigenvalues of $A$, this reduction method preserves stability. Furthermore, it preserves observability.

\begin{prop}[{\cite[Proposition 5]{TsuzukiOhki2024}}]\label{prop:inheritation_observability_controllability}
If the pair $(A, C)$ is observable, then the reduced pair $(A_{\bar{U}}, C_{\bar{U}})$ is also observable.
\end{prop}

By duality, if the pair $(A, B)$ is controllable, then projecting onto the dominant left invariant subspace preserves this property. That is, for $\bar{V} \in \mathcal{U}_r(A^\top)$, the reduced pair $(A_{\bar{V}}, B_{\bar{V}})$ is controllable. The inheritance of these properties thus depends on whether the projection is onto the right or left dominant invariant subspace.

This outcome can be verified using Gramians. 
Projecting the observability Gramian integral using $\bar{U}$ and the controllability Gramian using $\bar{V}$ for $T>0$ yields the Gramians of the respective reduced-order systems:
\begin{align}
	G_o(\bar{U},T) :=& \bar{U}^{\top} \left( \int_{0}^{T} \mathrm{e}^{A^{\top}t} C^{\top}C \mathrm{e}^{At} dt \right) \bar{U} 
	\nonumber \\
	=& \int_{0}^{T} \mathrm{e}^{A_{\bar{U}}^{\top}t} C_{\bar{U}}^{\top}C_{\bar{U}} \mathrm{e}^{A_{\bar{U}}t} dt, \label{eq:reduced_observability_grammian} 
	\\
	G_c(\bar{V},T) :=& \bar{V}^{\top} \left( \int_{0}^{T} \mathrm{e}^{At} BB^{\top} \mathrm{e}^{A^{\top}t} dt \right) \bar{V} 
	\nonumber \\ 
	=& \int_{0}^{T} \mathrm{e}^{A_{\bar{V}}t} B_{\bar{V}}B_{\bar{V}}^{\top} \mathrm{e}^{A_{\bar{V}}^{\top}t} dt. \label{eq:reduced_controllability_grammian}
\end{align}
Positive definiteness of the original Gramians implies positive definiteness of the projected ones, confirming the preservation of observability and controllability. However, these two reduced models, based on $\bar{U}$ and $\bar{V}$, exist in different state-space coordinate systems. 
To create a single reduced model that preserves both properties, we must relate them. 
Using the similarity transformation from Lemma \ref{lem:similar_transformation}, 
\begin{align*}
	& (\bar{V}^{\top} \bar{U} )^{-1} G_{c}(\bar{V},T) (\bar{V}^{\top} \bar{U} )^{-\top}
	\\
	= &\int_{0}^{T} \mathrm{e}^{A_{\bar{U}} t} \left( (\bar{V}^{\top} \bar{U} )^{-1} B_{\bar{V}} \right) \left( (\bar{V}^{\top} \bar{U} )^{-1} B_{\bar{V}} \right)^{\top} \mathrm{e}^{A_{\bar{U}}^{\top}t} dt
	.
	\end{align*}
	Because $G_{c}(\bar{V},T) >0$ implies $(\bar{V}^{\top} \bar{U} )^{-1} G_{c}(\bar{V},T) (\bar{V}^{\top} \bar{U} )^{-\top} >0$, the pair $( A_{\bar{U}}, B_{\bar{V}(\bar{V}^{\top}\bar{U})^{-1}})$ is controllable.
	
Similarly, similarity transformation of $G_{o}(\bar{U},T)$ by $(\bar{V}^{\top} \bar{U} )^{-\top}$ shows that if $G_{o}(\bar{U},T) >0$, the pair $( A_{\bar{V}}, C_{\bar{U}(\bar{V}^{\top}\bar{U})^{-1}})$ is observable. 
This argument is concluded in the following result:

\begin{prop}\label{prop:minimal_realization}
If the original system $(A,B,C)$ is a minimal realization (i.e., controllable and observable), then the following two reduced-order models are also minimal:
\begin{enumerate}
    \item $(A_{\bar{U}}, (\bar{V}^{\top}\bar{U})^{-1}B_{\bar{V}}, C_{\bar{U}})$
    \item $(A_{\bar{V}}, B_{\bar{V}}, C_{\bar{U}}(\bar{V}^{\top}\bar{U})^{-1})$
\end{enumerate}
\end{prop}

The two minimal models from Proposition \ref{prop:minimal_realization} are simply different state-space representations of the same input--output system.

\begin{prop}\label{prop:uniqueness_transfer_function}
The two reduced-order models from Proposition \ref{prop:minimal_realization} have the same transfer function.
\end{prop}

\begin{proof}
Let $R := \bar{V}^\top \bar{U}$. The transfer function for model 1 is $C_{\bar{U}}(sI_r - A_{\bar{U}})^{-1}R^{-1}B_{\bar{V}}$. Using the identity $A_{\bar{V}}R = R A_{\bar{U}}$ from Lemma \ref{lem:similar_transformation}, we have $(sI_r - A_{\bar{U}})^{-1}R^{-1} = R^{-1}(sI_r - A_{\bar{V}})^{-1}$. Substituting this yields $C_{\bar{U}}R^{-1}(sI_r - A_{\bar{V}})^{-1}B_{\bar{V}}$, which is the transfer function for model 2.
\end{proof}

Let $P_{\bar{U}}(s)$, $P_{\bar{V}}(s)$, and $P_{\rm rd}(s)$ be the transfer functions of the observability-preserving, controllability-preserving, and minimal reduced models, respectively. They are related by
\begin{align*}
	P_{\rm rd}(s) =& P_{\bar{U}}(s) + C_{\bar{U}} (sI_r - A_{\bar{U}})^{-1} (\bar{V}^\top \bar{U})^{-1} \bar{V}^\top \bar{U}_\perp B_{\bar{U}_\perp}
	\\
	=&
	P_{\bar{V}}(s) 
	+ 
	C_{\bar{V}_\perp} \bar{V}_\perp^{\top} \bar{U} (\bar{V}^{\top}\bar{U})^{-1} (sI_r - A_{\bar{V}})^{-1}  B_{\bar{V}}.
\end{align*}
This result shows that the minimal model $P_{\rm rd}(s)$ augments the observability-preserving model $P_{\bar{U}}(s)$ with a term that incorporates controllability information from the truncated subspace. A similar expression relates $P_{\rm rd}(s)$ to $P_{\bar{V}}(s)$. These three reduced models are generally distinct, as illustrated in the following example.

\begin{figure}[!htbp]
    	\centering
    	\includegraphics[keepaspectratio,width=\linewidth]{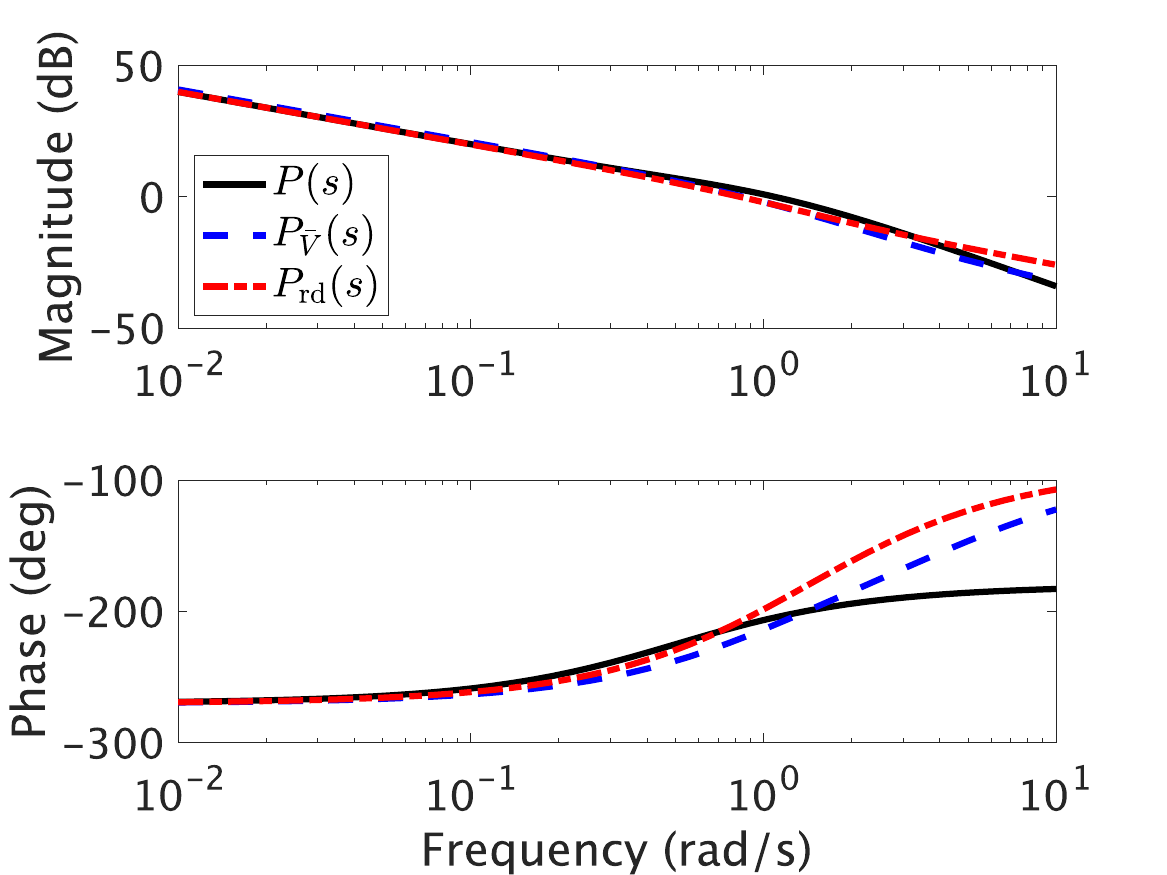}
    	\caption{Bode diagrams of $P(s)$, $P_{\bar{V}}(s)$, and $P_{\rm rd}(s)$. 
	}
	\label{fig:bode_plots}
	\end{figure}

\begin{exmp}\label{exmp:reduced_transfer_functions}
	Recall $A\in\mathbb{R}^{3\times 3}$ in Example \ref{exmp:St13} with the following matrices. 
	\begin{align*}
	A = \begin{bmatrix} 1 & 1 & 2 \\ 0 & 0 & 1 \\ 0 & 0 & -1 \end{bmatrix}
	,\
	B = \begin{bmatrix} 0 \\ 0 \\ 1 \end{bmatrix} ,\ 
	C = \begin{bmatrix} 1 & 0 & 0 \end{bmatrix}
	.
	\end{align*}
	The transfer function $P(s)$ of the original system with $(A,B,C)$ and the reduced transfer functions of the model $(A_{\bar{U}}, B_{\bar{U}},C_{\bar{U}})$ and $(A_{\bar{V}}, B_{\bar{V}},C_{\bar{V}})$ for $\bar{U} \in \mathcal{U}_{2}(A) $ and $\bar{V} \in \mathcal{U}_{2}(A^{\top})$ are 
	\begin{align*}	
	P(s) =& \frac{2s+1}{s^{3}-s}
	,\quad 
	P_{\bar{U}} (s) =0, \quad 
	P_{\bar{V}} (s) = \frac{1}{9}\frac{2s+1}{s (s-1)}
	.
	\end{align*}	
	Note that 
	\begin{align*}
	\psi_{1}(A^{\top}) = \frac{1}{\sqrt{17}}\begin{bmatrix}
	2 \\ 2 \\ 3
	\end{bmatrix}
	, 
	\psi_{2}(A^{\top}) = \frac{1}{\sqrt{2}}\begin{bmatrix}
	0 \\ 1 \\ 1
	\end{bmatrix}
	, 
	\psi_{3}(A^{\top}) =\begin{bmatrix}
	0 \\ 0 \\ 1
	\end{bmatrix}
	.
	\end{align*}
	Because the eigenvectors $\psi _{1}(A)$ and $\psi_{2}(A)$ denoted in Example \ref{exmp:St13} are orthogonal to $B$, $P_{\bar{U}}(s) =0$ and $P_{\bar{V}}(s) \neq 0$. 
	
	Conversely, the reduced transfer function of the model $(A_{\bar{U}}, B_{\bar{V}(\bar{V}^{\top} \bar{U})^{-\top}} , C_{\bar{U}} )$ is 
	\begin{align*}
	P_{\rm rd} (s) = \frac{1}{2} \frac{2s+1}{s(s-1)}.
	\end{align*}
	In this example, each reduced transfer function is different from the others. 
	Fig. \ref{fig:bode_plots} shows the Bode diagrams of $P(s)$, $P_{\bar{V}}(s)$, and $P_{\rm rd}(s)$.  
	$P_{\bar{V}}(s)$ and $P_{\rm rd}(s)$ well approximate $P(s)$ in the low-frequency range, and the phase differs in the high-frequency range. 
	A detailed analysis of the approximation error will be studied in the future. 
	\end{exmp}
	
	Although we have a controllability- and observability-preserving reduced model $(A_{\bar{U}}, B_{\bar{V}(\bar{V}^{\top} \bar{U})^{-\top}}, C_{\bar{U}} )$, as mentioned in Example \ref{exmp:reduced_transfer_functions}, it is unclear whether the reduction is more effective than the other ones. 
	For example, as demonstrated here, the others have intuitive approximation errors: 
	
	\begin{prop} \label{prop:approximation_error_freq}
	Consider a system with $(A,B,C)$.  
	Then, for $\bar{U} \in \mathcal{U}_{r}(A)$ and $\bar{V} \in \mathcal{U}_{r}(A^{\top})$, the transfer function is 
	\begin{align*}
	P(s) =& P_{\bar{U}}(s) 
	+
	C_{\bar{U}_{\perp}} (s I_{n-r} - A_{\bar{U}_{\perp}})^{-1} B_{\bar{U}_{\perp}}
	\\
	&
	+ C_{\bar{U}} (sI_{r} - A_{\bar{U}})^{-1} \bar{U}^{\top} A \bar{U}_{\perp} ( sI_{n-r} -A _{\bar{U}_{\perp}} ) ^{-1} B_{\bar{U}_{\perp}}
	\\
	=& P_{\bar{V}}(s) 
	+
	C_{\bar{V}_{\perp}} (s I_{n-r} - A_{\bar{V}_{\perp}})^{-1} B_{\bar{V}_{\perp}}
	\\
	&
	+ C_{\bar{V}_{\perp}} ( sI_{n-r} -A _{\bar{V}_{\perp}} )^{-1}  \bar{V}_{\perp}^{\top} A \bar{V}(sI_{r} - A_{\bar{U}})^{-1} B_{\bar{V}}
	,
	\end{align*}
	where $\bar{U}_{\perp} \in \mathcal{U}_{n-r} (-A^{\top})$ and $\bar{V} _{\perp} \in \mathcal{U}_{n-r} (-A)$. 
	\end{prop}
	
	Recall that $\bar{U}_{\perp} \bar{U}_{\perp} ^{\top} = I_{n} - \bar{U} \bar{U} $ and $\bar{V} _{\perp} \bar{V}_{\perp}^{\top} = I_{n} - \bar{V} \bar{V}^{\top}$. 
	
	\begin{proof}
	Using \eqref{eq:blk_upper_triangle}, 
	\begin{align*}
	&P(s) 
	\\
	=& 
	C Q_{\bar{U}} (sI_{n} - Q_{\bar{U}} ^{\top}AQ_{\bar{U}} ) ^{-1} Q_{\bar{U}} ^{\top} B
	\\
	=&
	\begin{bmatrix}
	C_{\bar{U}}^{\top} \\ C_{\bar{U}_{\perp}}^{\top}
	\end{bmatrix}^{\top}
	\begin{bmatrix}
	(sI_{r} - A_{\bar{U}})^{-1}  & P_{ur}(s) 
	\\
	0 _{n-r,r} & ( s I_{n-r} - A_{\bar{U}_{\perp}} )^{-1}
	\end{bmatrix} 
	\begin{bmatrix}
	B_{\bar{U}} \\ B_{\bar{U}_{\perp}}
	\end{bmatrix},
	\end{align*}
	where the upper right block is $P_{ur} (s) := (sI_{r} - A_{\bar{U}}) ^{-1} \bar{U}^{\top} A \bar{U}_{\perp} (s I_{n-r} - A_{\bar{U}_{\perp}} ) ^{-1}$. 
	Hence, 
	\begin{align*}
	P(s)=&
	C_{\bar{U}} (sI_{r} - A_{\bar{U}})^{-1} B_{\bar{U}}
	+
	C_{\bar{U}_{\perp}} (s I_{n-r} - A_{\bar{U}_{\perp}})^{-1} B_{\bar{U}_{\perp}}
	\\
	&
	+ C_{\bar{U}}(sI_{r} - A_{\bar{U}}) ^{-1} \bar{U}^{\top} A \bar{U}_{\perp} (s I_{n-r} - A_{\bar{U}_{\perp}} ) ^{-1} B_{\bar{U}_{\perp}}
	.
	\end{align*}
	Similarly, 
	\begin{align*}
	P(s) = &C_{\bar{V}} (sI_{r} - A_{\bar{V}})^{-1} B_{\bar{V}}
	+
	C_{\bar{V}_{\perp}} (s I_{n-r} - A_{\bar{V}_{\perp}})^{-1} B_{\bar{V}_{\perp}}
	\\
	&
	+ C_{\bar{V}_{\perp}} ( sI_{n-r} -A _{\bar{V}_{\perp}} ) ^{-1}  \bar{V}_{\perp}^{\top} A \bar{V}(sI_{r} - A_{\bar{V}})^{-1} B_{\bar{V}}
	\end{align*}
	holds. 
	\end{proof}
	
\begin{figure}[!htbp]
    	\centering
    	\includegraphics[keepaspectratio,width=\linewidth]{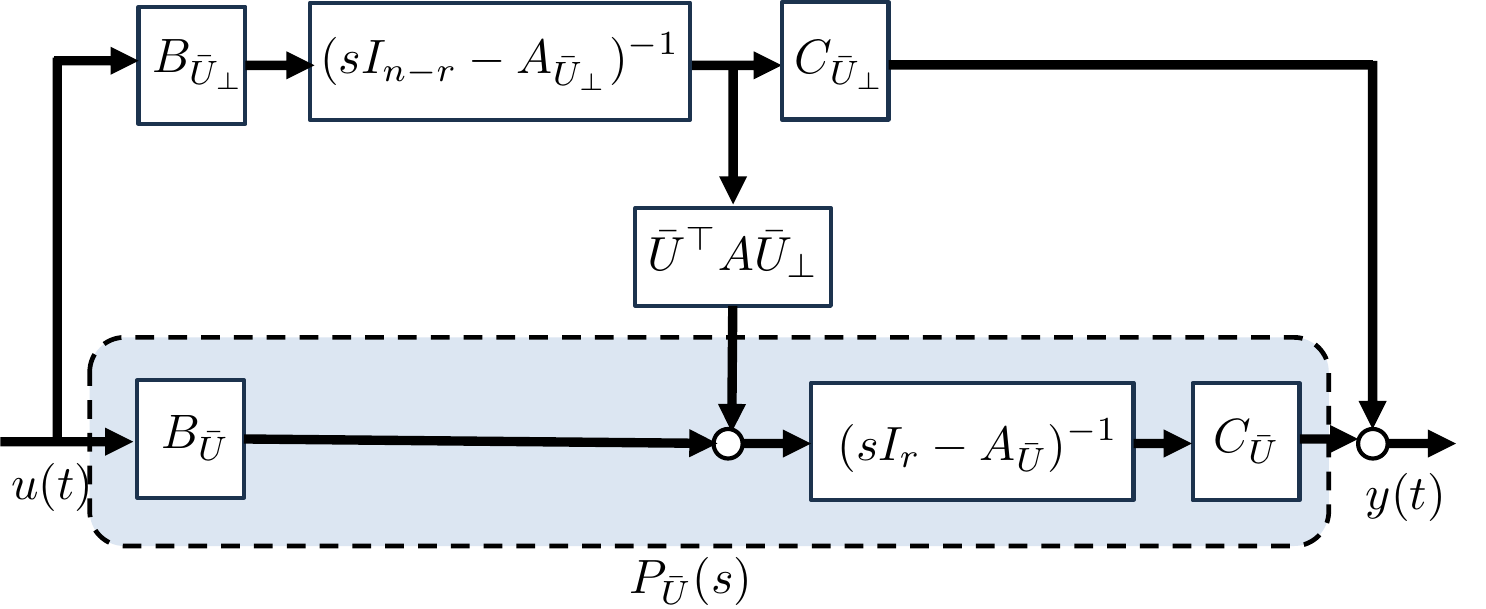}
    	\caption{Block diagram of the transfer function $P(s)$ and $P_{\bar{U}}(s)$. 
	}
	\label{fig:Oja_based_model_reduction}
	\end{figure}

	Fig. \ref{fig:Oja_based_model_reduction} shows the block diagram of $P(s)$ and $P_{\bar{U}}(s)$ in Proposition \ref{prop:approximation_error_freq}. 
	If one considers the approximation error $P(s) - P_{\rm rd}(s)$, it is also calculated from the relationship between $P_{\rm rd}(s)$ and $P_{\bar{U}}(s)$; we skip it. 
	The error includes a cascade-connected term from the truncated part to the reduced model, which means that the error dynamics is unstable if the original system is unstable.  
	If the original system is stable, the approximation error system is also stable.  
	Hence, we can establish the approximation error by the well-known $H^{2}$ norm or $H^{\infty}$ norm. 
	
	\begin{remark}[On the Minimality of Reduced Models] \label{rem:frequency_domain_error_analysys}
If the original system $(A,B,C)$ is a minimal realization, then the truncated input matrix $B_{\bar{U}_{\perp}}$ must be non-zero. If $B_{\bar{U}_{\perp}} = 0$, the error formula in Proposition \ref{prop:approximation_error_freq} would simplify to $P(s) = P_{\bar{U}}(s)$. This result would imply that the original $n$-th order system has a lower-order, $r$-th dimensional representation, which contradicts the initial assumption of minimality.

Furthermore, if the original system is minimal and the coupling term $\bar{U}^{\top} A \bar{U}_{\perp}$ is zero, then the reduced-order model with transfer function $P_{\bar{U}}(s)$ is also a minimal realization. When this coupling term is zero, the system decomposes perfectly, and the overall transfer function becomes a direct sum: $P(s) = P_{\bar{U}}(s) + P_{\bar{U}_{\perp}}(s)$. The McMillan degree (i.e., the order) of this sum is the sum of the degrees of its parts, assuming no pole-zero cancellations. If the realization for $P_{\bar{U}}(s)$ were not minimal, its degree would be less than $r$, causing the total degree of $P(s)$ to be less than $n$. This result again contradicts the minimality of the original system. By duality, the same holds for $P_{\bar{V}}(s)$ if $\bar{V}_{\perp}^{\top}A \bar{V} = 0$.
\end{remark}

	\subsection{Low-Rank Design for Observer and State Feedback Gains}

	Controller synthesis for large-scale systems often relies on model reduction. However, feedback control designed for a reduced-order model must be robust to the dynamics of the truncated part of the system, often necessitating complex robust control techniques such as $H^\infty$ synthesis. In this section, we present an alternative approach: synthesizing low-rank controllers that stabilize the full-order system directly, thereby avoiding truncation errors.

Consider the LTI system \eqref{eq:system}. We assume that the system matrix $A$ has at most $r$ unstable eigenvalues, where $r \ll n$, and that the eigenvalue gap condition $\mathrm{Re}(\lambda_r(A)) > \mathrm{Re}(\lambda_{r+1}(A))$ holds. Although full-order observer-based control is a fundamental stabilization strategy, its design and implementation become challenging for large $n$. The Oja flow-based method presented here offers a computationally tractable design process by focusing only on the low-dimensional unstable subspace.

From the results discussed in the previous section, the reduced-order pairs $(A_{\bar{U}}, C_{\bar{U}})$ and $(A_{\bar{V}}, B_{\bar{V}})$ inherit the observability and controllability of the original system's unstable modes. Therefore, if the unstable modes of $(A,B,C)$ are stabilizable and detectable, there exist gain matrices $L_r \in \mathbb{R}^{r \times p}$ and $F_r \in \mathbb{R}^{m \times r}$ such that the $r \times r$ matrices
\begin{align*}
	A_{\bar{U}} - L_r C_{\bar{U}} \quad \text{and} \quad A_{\bar{V}} - B_{\bar{V}} F_r
\end{align*}
are Hurwitz. The key insight is that these low-dimensional gains can be embedded back into the full-dimensional space to stabilize the original system.

\begin{prop}[Low-Rank Stabilization] \label{prop:low_dimensional_FB_design}
If the $r \times r$ matrices $A_{\bar{U}} - L_r C_{\bar{U}}$ and $A_{\bar{V}} - B_{\bar{V}} F_r$ are Hurwitz, then the full-order, low-rank closed-loop matrices
\begin{align*}
	A - \bar{U}L_r C \quad \text{and} \quad A - BF_r \bar{V}^\top
\end{align*}
are also Hurwitz.
\end{prop}

\begin{proof}
	The proof follows the arguments in \cite[Proposition 8]{TsuzukiOhki2024}. Consider the observer error dynamics matrix $A - \bar{U}L_r C$. 
	As shown in \eqref{eq:blk_upper_triangle}, a similarity transformation by the orthogonal matrix $Q_{\bar{U}} = [\bar{U}, \bar{U}_\perp]$ yields a block upper-triangular matrix:
\begin{align*}
	Q_{\bar{U}}^{\top}(A - \bar{U}L_r C) Q_{\bar{U}} = \begin{bmatrix}
	A_{\bar{U}} - L_r C_{\bar{U}} & (\bar{U}^{\top} A - L_r C)\bar{U}_{\perp} \\
	0_{n-r,r} & A_{\bar{U}_\perp}
	\end{bmatrix}.
\end{align*}
The eigenvalues of this matrix are the union of the eigenvalues of the diagonal blocks. The top-left block, $A_{\bar{U}} - L_r C_{\bar{U}}$, is Hurwitz by design. The bottom-right block, $A_{\bar{U}_\perp}$, contains the stable eigenvalues $\{\lambda_{r+1}(A), \dots, \lambda_n(A)\}$. Because all eigenvalues have negative real parts, the full matrix $A - \bar{U}L_r C$ is Hurwitz.
The stability of $A - BF_r\bar{V}^\top$ follows from a symmetric argument by applying the same logic to $(A^\top, B^\top, \bar{V}, F_r^\top)$.
\end{proof}

Proposition \ref{prop:low_dimensional_FB_design} enables the design of a stabilizing observer-based controller. A Luenberger observer for the system is given by
\begin{align*}
	\frac{d}{dt}\hat{x}(t) = A\hat{x}(t) + Bu(t) + \bar{U}L_r(y(t) - C\hat{x}(t)),
\end{align*}
where $\hat{x}(t) \in \mathbb{R}^n$. The estimation error $e(t) := x(t) - \hat{x}(t)$ evolves according to $\frac{d}{dt}e(t) = (A - \bar{U}L_r C)e(t)$. Because this system matrix is Hurwitz, the error converges to zero exponentially.

Applying the state-feedback control law $u(t) = -F_r\bar{V}^\top \hat{x}(t)$ yields the full closed-loop dynamics for the state $x(t)$ and error $e(t)$:
\begin{align*}
	\frac{d}{dt} \begin{bmatrix} x(t) \\ e(t) \end{bmatrix} =
	\begin{bmatrix}
	A - BF_r\bar{V}^\top & BF_r\bar{V}^\top \\
	0_{n} & A - \bar{U}L_r C
	\end{bmatrix}
	\begin{bmatrix} x(t) \\ e(t) \end{bmatrix}.
\end{align*}
Because of the block-triangular structure, the stability of the overall system is guaranteed by the stability of the diagonal blocks, which are both Hurwitz by Proposition \ref{prop:low_dimensional_FB_design}.

Although the implementation of this observer is still $n$-dimensional, the computationally intensive part of the design process—--finding the gains $L_r$ and $F_r$ by pole placement or solving Riccati equations—--is performed on the small $r \times r$ systems. 
The low-rank approximated Kalman--Bucy filter \cite{TsuzukiOhki2024} is an example of this result.  
This result renders the design of stabilizing controllers for very large-scale systems with few unstable modes computationally feasible.

	\subsection{Singularly Perturbed Linear Systems}
	
	Consider a stable LTI system whose dynamics exhibit a two-time-scale property. This system is characterized by a Hurwitz matrix $A$ whose eigenvalues are separated into a group of $r$ {\em slow} and $n-r$ {\em fast} eigenvalues, such that there is a significant gap between them: $|\mathrm{Re}(\lambda_r(A))| \ll |\mathrm{Re}(\lambda_{r+1}(A))|$. The Oja flow provides a natural method for separating these modes. The subspace spanned by $\bar{U} \in \mathcal{U}_r(A)$ corresponds to the slow dynamics, whereas the orthogonal complement, spanned by $\bar{U}_\perp$, corresponds to the fast dynamics.

This separation is the foundation of the singular perturbation method, which decouples a system to derive a simplified model of the slow dynamics \cite{Kokotovic_optimal_control}. Applying this idea to system \eqref{eq:system}, we define the slow states as $x_s(t) := \bar{U}^\top x(t)$ and the fast states as $x_f(t) := \bar{U}_\perp^\top x(t)$. With the block-triangular decomposition from \eqref{eq:blk_upper_triangle}, the state equations become
\begin{align*}
	\frac{d}{dt} x_s(t) &= A_{\bar{U}} x_s(t) + (\bar{U}^\top A \bar{U}_\perp) x_f(t) + B_{\bar{U}} u(t), \\
	\frac{d}{dt} x_f(t) &= A_{\bar{U}_\perp} x_f(t) + B_{\bar{U}_\perp} u(t).
\end{align*}
The core assumption of singular perturbation is that the fast dynamics reach a quasi-steady state much more rapidly than the slow dynamics evolve, provided the input $u(t)$ also varies slowly. By setting $\frac{d}{dt}x_f(t) \approx 0$, we can solve for the quasi-steady-state value of the fast mode: $x_f(t) \approx -A_{\bar{U}_\perp}^{-1} B_{\bar{U}_\perp} u(t)$. Substituting this back into the equation for the slow dynamics yields the reduced-order model:
\begin{align*}
	\frac{d}{dt} x_s(t) = A_{\bar{U}} x_s(t) + \left( B_{\bar{U}} - (\bar{U}^\top A \bar{U}_\perp) A_{\bar{U}_\perp}^{-1} B_{\bar{U}_\perp} \right) u(t).
\end{align*}
Although the singular perturbation method is powerful, it can be challenging to identify the transformation that separates the slow and fast modes for large and complex systems. The Oja flow offers a significant advantage by providing a systematic and computationally tractable way to find this decomposition automatically.

\section{Conclusion}
	
	This paper presents a global convergence analysis of Oja's component flow for extracting the $r$-dominant invariant subspace of a general square matrix $A$, which corresponds to the $r$ eigenvalues with the largest real parts. The stable invariant set, $\mathcal{U}_r$, and its domain of attraction were characterized. Crucially, the algorithm was shown to ensure exponential convergence for almost all initial conditions, providing a strong theoretical guarantee for its reliability. We also investigated the tracking performance for time-varying matrices and provided numerical examples validating the theoretical results.

The primary implication of these findings is that the Oja flow is a theoretically sound and practically viable tool for analyzing dominant dynamics in a much broader class of linear systems than previously considered, extending beyond symmetric positive-definite matrices. The applications presented demonstrate this utility. We propose methods for model reduction and low-rank controller synthesis that leverage the Oja flow's ability to isolate key subspaces. The properties of the resulting reduced-order models, including preservation of observability, controllability, and stabilization using output feedback are discussed. Furthermore, we demonstrated the flow's applicability to mode identification in singularly perturbed systems by exploiting its capacity to separate slow and fast dynamics. These concepts offer promising extensions to time-varying and discrete-time systems.

Future research directions include: 1) accelerating the convergence rate for enhanced practical performance; 2) establishing rigorous analytical guarantees for time-varying matrices, essential for model reduction of more general systems; and 3) further characterizing the properties of the derived reduced-order models. Regarding the time-varying case, although applications to state estimation exist \cite{tranninger2022detectability}, performance is sensitive to the parameter $\varepsilon$, often requiring small integration steps and increasing computational cost. A more detailed analysis of tracking error is therefore required, alongside broader exploration of estimation and control applications. Finally, refining the understanding of the reduced-model properties remains a challenging yet important problem because of its direct impact on practical controller design and performance guarantees.

%

\end{document}